\documentclass[12pt,final]{elsart}
\usepackage{amsmath,amssymb,amsfonts}
\usepackage[utf8]{inputenc} 
\usepackage[T1]{fontenc}
\usepackage{enumerate}
\usepackage{pdfpages}
\usepackage{mathrsfs}
\usepackage{graphicx}
\usepackage{tabularx}
\usepackage{xcolor}
\usepackage{url}
\newcommand{\MyBibitem}[2]{\bibitem{#1}}
\newcommand{\MyReference}[3]{\textsc{#1}. \emph{#2}. #3.}
\newcommand{\MyParagraph}[1]{\textbf{#1.}\;}

\newcommand{\NN}{\mathbb{N}}
\newcommand{\ZZ}{\mathbb{Z}}
\newcommand{\RR}{\mathbb{R}}
\newcommand{\CC}{\mathbb{C}}
\newcommand{\dd}{\text{d}}
\newcommand{\ee}{\text{e}}
\newcommand{\ii}{\text{i}}
\newcommand{\abs}[1]{|#1|}

\newcommand{\nE}{\mathcal{E}}
\newcommand{\nF}{\mathcal{F}}
\newcommand{\nM}{\mathcal{M}}
\newcommand{\nO}{\mathcal{O}}
\begin{document}
\begin{frontmatter}
\title{Generalisation of splitting methods based \\ on modified potentials to nonlinear evolution equations of parabolic and Schr{\"o}dinger type}
\author{S.~Blanes, F.~Casas, C.~Gonz{\'a}lez, M.~Thalhammer}%
\footnote{\scriptsize $\!\!\!\!$Addresses:
Sergio Blanes, Universitat Polit{\`e}cnica de Val{\`e}ncia, Instituto de Matem{\'a}tica Multi\-disciplinar, 46022~Valencia, Spain. 
Fernando Casas, Universitat Jaume~I, IMAC and Departament de Matem{\`a}tiques, 12071~Castell\'on, Spain. 
Ces{\'a}reo Gonz{\'a}lez, Universidad de Valladolid, Departamento de Matem{\'a}tica Aplicada, 47011~Valladolid, Spain. 
Mechthild Thalhammer, Leopold--Franzens-Universit{\"a}t Innsbruck, Institut f{\"u}r Mathematik, 6020~Innsbruck, Austria.
Websites: \url{personales.upv.es/serblaza}, \url{www.gicas.uji.es/fernando.html}, 
\url{www.imuva.uva.es/en/investigadores/33}, \url{techmath.uibk.ac.at/mecht}. 
Email addresses: \url{serblaza@imm.upv.es}, \url{fernando.casas@uji.es}, \url{ome@am.uva.es}, \url{mechthild.thalhammer@uibk.ac.at}.}
\address{}
\begin{abstract}
{The present work is concerned with the extension of modified potential operator splitting methods to specific classes of nonlinear evolution equations.
The considered partial differential equations of Schr{\"o}dinger and parabolic type comprise the Laplacian, a potential acting as multiplication operator, and a cubic nonlinearity. 
Moreover, an invariance principle is deduced that has a significant impact on the efficient realisation of the resulting modified operator splitting methods for the Schr{\"o}dinger case.}

Numerical illustrations for the time-dependent Gross--Pitaevskii equation in the physically most relevant case of three space dimensions and for its parabolic counterpart related to ground state and excited state computations confirm the benefits of the proposed fourth-order modified operator splitting method in comparison with standard splitting methods.

The presented results are novel and of particular interest from both, a theoretical perspective to inspire future investigations of modified operator splitting methods for other classes of nonlinear evolution equations and a practical perspective to advance the reliable and efficient simulation of Gross--Pitaevskii systems in real and imaginary time.
\end{abstract}
\end{frontmatter}
\section{Introduction}
\label{sec:Section1}
\MyParagraph{Scope of applications}
A wide range of relevant applications in sciences includes the numerical integration of initial value problems for nonlinear evolution equations.
In many cases, the function defining the right-hand side comprises two or more parts
\begin{equation}
\label{eq:IVPNonlinear}
\begin{cases}
\tfrac{\dd}{\dd t} \, u(t) = F_1\big(u(t)\big) + F_2\big(u(t)\big)\,, \\
u(0) = u_0\,, \quad t \in [0, T]\,.
\end{cases}
\end{equation}
As prominent instances, we highlight nonlinear Schr{\"o}dinger equations, more specifically, time-dependent Gross--Pitaevskii equations that arise in the description of Bose--Einstein condensation, see~\cite{Gross1961,Pitaevskii1961}.
For a comprehensive overview of the underlying principles of quantum theory, we refer to~\cite{Messiah1999}.

\MyParagraph{Nonlinear Schr{\"o}dinger equation (Gross--Pitaevskii equation)}
A fundamental model for the nonlinear dynamics of a single Bose--Einstein condensate reads as 
\begin{subequations}
\label{eq:GPE}  
\begin{equation}
\label{eq:GPE1}  
\begin{cases}
\ii \, \partial_t \Psi(x, t) = - \, \Delta \Psi(x, t) + V(x) \, \Psi(x, t) + \vartheta \, \abs{\Psi(x, t)}^2 \, \Psi(x, t)\,, \\
\Psi(x, 0) = \Psi_0(x)\,, \quad (x, t) \in \Omega \times [0, T]\,, 
\end{cases}
\end{equation}  
where $\Delta = \partial_{x_1}^2 + \dots + \partial_{x_d}^2$ denotes the Laplacian with respect to the spatial variables $x = (x_1, \dots, x_d) \in \RR^d$, $V: \RR^d \to \RR$ a real-valued potential, $\vartheta \in \RR$ the coupling constant, and $\Psi: \Omega \times [0, T] \subset \RR^d \times \RR \to \CC$ the space-time-dependent complex-valued macroscopic wave function. 
Assigning for a regular function $v: \Omega \to \CC$ the linear differential and nonlinear multiplication operators
\begin{equation}
\label{eq:F1F2}
\begin{gathered}
\big(F_1(v)\big)(x) = c \, \Delta \, v(x)\,, \quad c = \ii\,, \\
\big(F_2(v)\big)(x) = \bar{c} \, \big(V(x) + \vartheta \, \abs{v(x)}^2\big) \, v(x)\,, \quad \bar{c} = - \, \ii\,, \\
x \in \Omega\,, 
\end{gathered}
\end{equation}
\end{subequations}
and setting $u(t) = \Psi(\cdot, t)$ for $t \in [0, T]$, we retain the general formulation~\eqref{eq:IVPNonlinear}.

\MyParagraph{Nonlinear parabolic equation}
By analogy to the time-dependent Gross--Pitaevskii equation~\eqref{eq:GPE1}, we consider the parabolic problem 
\begin{subequations}
\label{eq:Parabolic}    
\begin{equation}
\label{eq:Parabolic1}    
\begin{cases}
\partial_t U(x, t) = \Delta U(x, t) + V(x) \, U(x, t) + \vartheta \, \abs{U(x, t)}^2 \, U(x, t)\,, \\
U(x, 0) = U_0(x)\,, \quad (x, t) \in \Omega \times [0, T]\,,
\end{cases}
\end{equation}  
for a real-valued solution $U: \Omega \times [0, T] \subset \RR^d \times \RR \to \RR$.
Accordingly, it corresponds to~\eqref{eq:F1F2} with different constant
\begin{equation}
\begin{gathered}
\big(F_1(v)\big)(x) = c \, \Delta \, v(x)\,, \quad c = 1\,, \\
\big(F_2(v)\big)(x) = \bar{c} \, \big(V(x) + \vartheta \, \abs{v(x)}^2\big) \, v(x)\,, \quad \bar{c} = 1\,, \\
x \in \Omega\,, 
\end{gathered}
\end{equation}
\end{subequations}
and, setting $u(t) = U(\cdot, t)$ for $t \in [0, T]$, we obtain again the general form~\eqref{eq:IVPNonlinear}.
It is noteworthy that the parabolic equation~\eqref{eq:Parabolic1} arises in ground state and excited state computations, see for instance~\cite{BaoDu2004,DanailaProtas2017}.

{\MyParagraph{Splitting approach}
In essence, operator splitting methods rely on the presumption that the numerical approximation of the subproblems 
\begin{equation*}
\begin{gathered}
\tfrac{\dd}{\dd t} \, u_1(t) = F_1\big(u_1(t)\big)\,, \quad
\tfrac{\dd}{\dd t} \, u_2(t) = F_2\big(u_2(t)\big)\,,
\end{gathered}
\end{equation*}  
is significantly simpler compared to the numerical approximation of the original problem~\eqref{eq:IVPNonlinear}.
Then, within multiple scopes,}
for ordinary differential equations and time-dependent partial differential equations, for linear problems as well as nonlinear problems, a variety of works has confirmed the benefits of operator splitting methods regarding desirable features that are subsumed under the central concepts stability, efficiency, and preservation of conserved quantities.
For general information, we refer to~\cite{McLachlanQuispel2002,SanzSernaCalvo2018}.
Specific studies in the context of Schr{\"o}dinger equations are given, e.g., in~\cite{BaoJinMarkowich2002,BlanesCasasMurua2015,CaliariZuccher2021,ThalhammerAbhau2012}.

{\MyParagraph{Alternative approach}
In this work, we propose an approach that provides a favourable alternative to standard operator splitting methods in situations, where the operator~$F_2$ and an iterated commutator of~$F_2$ and~$F_1$, given by 
\begin{equation}
\label{eq:G2Introduction}
\begin{split}
G_2(v) &= F_1''(v) \, F_2(v) \, F_2(v) + F_1'(v) \, F_2'(v) \, F_2(v) + F_2'(v) \, F_2'(v) \, F_1(v) \\
&\qquad - F_2''(v) \, F_1(v) \, F_2(v) - 2 \, F_2'(v) \, F_1'(v) \, F_2(v)\,, 
\end{split}
\end{equation}  
have a similar structure.
As relevant nonlinear partial differential equations with this property, we identify Schr{\"o}dinger and parabolic equations such as~\eqref{eq:GPE} and~\eqref{eq:Parabolic} that comprise the Laplacian, a potential acting as multiplication operator, and a cubic nonlinearity. 
For specifications concerning~\eqref{eq:G2Introduction}, we in particular refer to Sections~\ref{sec:Section3} and~\ref{sec:Section4}. 
}

{\MyParagraph{Formal means and objectives}
Our educated guess that leads us to modified operator splitting methods relies on a formal generalisation of the linear case, which we briefly sketch next and describe in further detail in the subsequent sections. 
\begin{enumerate}[(i)]
\item
\emph{Linear ordinary differential equations.} \;   
The starting point is a linear ordinary differential equation defined by non-commuting square matrices
\begin{equation*}
\tfrac{\dd}{\dd t} \, u(t) = A \, u(t) + B \, u(t)\,, \quad t \in [0, T]\,.
\end{equation*}
The corresponding solution value at the final time is given by the matrix exponential, that is   
\begin{equation*}
u(T) = \ee^{T (A + B)} \, u(0) = \Big(\ee^{\tau (A + B)}\Big)^N \, u(0)\,, \quad \tau = \tfrac{T}{N}\,, \quad N \in \NN\,.
\end{equation*}
Standard splitting methods are built on compositions of the factors~$\ee^{a \, \tau A}$ and~$\ee^{b \, \tau B}$ with suitably chosen coefficients~$a$ and~$b$.
Beyond that, components of the form
\begin{equation}
\label{eq:BBAIntroduction}
\ee^{b \, \tau B + c \, \tau^3 \, [B,[B,A]]}\,, \quad \big[B,[B,A]\big] = B^2 A - 2 \, B A B + A B^2\,, 
\end{equation}
with certain coefficients~$b$ and~$c$ are incorporated in modified potential operator splitting methods.
The underlying idea of this approach is to gain freedom in the adjustment of the method coefficients and, amongst others,
to overcome an order barrier valid for standard splitting methods. 
\item
\emph{Linear partial differential equations.} \;   
Advantages of this approach become apparent in the context of the imaginary time integration of linear Schr{\"o}dinger equations comprising the Laplacian and a potential.
There, the operator arising in~\eqref{eq:BBAIntroduction} reduces to a multiplication operator, which is defined by the potential and its gradient. 
\item
\emph{Nonlinear partial differential equations.} \;   
The guide line for the extension to nonlinear evolution equations~\eqref{eq:IVPNonlinear} is provided by the calculus of Lie derivatives, see~\cite{Varadarajan1984} for a detailed exposition.
In order to make our contribution accessible to a broader readership, we do not presume the knowledge of this formal calculus and explain the required elementary means on occasion.
Basically, the operators~$F_1$ and~$F_2$ take the roles of the matrices~$A$ and~$B$, and~\eqref{eq:BBAIntroduction} is replaced by the solution to 
\begin{equation*}
\tfrac{\dd}{\dd t} \, u(t) = b \, F_2\big(u(t)\big) + c \, \tau^2 \, G_2\big(u(t)\big)\,, \quad t \in [0, \tau]\,, 
\end{equation*}  
see also~\eqref{eq:G2Introduction}.
This formalism is expedient with regard to the design of novel higher-order time integration methods for nonlinear partial differential equations. 
Though, it is then equally of importance to concretise formal considerations and to confirm that the resulting modified operator splitting methods are indeed practicable and beneficial. 
\item
\emph{Main objectives.} \;   
In this work, for the sake of concreteness, we focus on the extension of a famous fourth-order modified potential operator splitting method by \textsc{Chin}~\cite{Chin1997} to the Gross--Pitaevskii equation~\eqref{eq:GPE} and the parabolic equation~\eqref{eq:Parabolic}.
So far, this scheme has been introduced and studied merely for the linear case.
\end{enumerate}
}

{\MyParagraph{Outline}
The present manuscript is organised as follows. 
In Section~\ref{sec:Section2}, we review fundamental concepts for operator splitting methods.
In Sections~\ref{sec:Section3} and~\ref{sec:Section4}, we state the formal generalisation of a fourth-order modified potential operator splitting method to the nonlinear case and substantiate it for the Gross--Pitaevskii equation~\eqref{eq:GPE} and its parabolic analogue~\eqref{eq:Parabolic}.
A fundamental invariance principle that includes a known result for standard splitting methods as a special case is deduced in Section~\ref{sec:Section5}. 
In Section~\ref{sec:Section6}, we detail the implementation of the novel modified operator splitting method based on a Fourier spectral space discretisation and provide numerical comparisons with standard splitting methods.
Additional information on a publicly accessible \textsc{Matlab} code is found in Appendix~\ref{sec:AppendixMatlabCode}. 
The observed order reduction of Yoshida's fourth-order complex splitting method~\eqref{eq:CoefficientsYoshidaComplex} is analysed in Appendix~\ref{sec:AppendixOrderReduction}.
}
\section{{Survey of standard and modified potential splitting methods}}
\label{sec:Section2}
\MyParagraph{Linear case}
As an illustrative example, we state the simplest representative of standard splitting methods, the first-order Lie--Trotter splitting method, for a system of linear differential equations 
\begin{equation}
\label{eq:IVPLinear}
\begin{cases}
\tfrac{\dd}{\dd t} \, u(t) = A \, u(t) + B \, u(t)\,, \\
u(0) = u_0\,, \quad t \in [0, T]\,, 
\end{cases}
\end{equation}
defined by non-commuting time-independent complex matrices $A, B \in \CC^{M \times M}$.
For a positive integer number $N \in \NN$ with associated time increment and equidistant grid points
\begin{equation*}
\tau = \tfrac{T}{N}\,, \quad t_n = n \, \tau\,, \quad n \in \{0, 1, \dots, N\}\,, 
\end{equation*}
numerical approximations to the exact solution values are obtained by the recurrence
\begin{equation*}
u_{n+1} = \ee^{\tau B} \, \ee^{\tau A} \, u_n \approx u(t_{n+1})\,, \quad n \in \{0, 1, \dots, N-1\}\,.
\end{equation*}
Higher-order splitting methods for~\eqref{eq:IVPLinear} involve the action of several matrix exponentials on the current approximation and can be cast into the format 
\begin{equation*}
u_{n+1} = \ee^{b_s \tau B} \, \ee^{a_s \tau A} \, \cdots \, \ee^{b_1 \tau B} \, \ee^{a_1 \tau A} \, u_n \approx u(t_{n+1})\,, \quad n \in \{0, 1, \dots, N-1\}\,,
\end{equation*}
with real or complex coefficients $(a_j, b_j)_{j=1}^s$, respectively.

\MyParagraph{Nonlinear case}
Their generalisation to nonlinear evolution equations~\eqref{eq:IVPNonlinear} is based on the composition of the solutions to the subproblems defined by~$F_1$ and~$F_2$.
We henceforth employ the compact notation
\begin{subequations}
\label{eq:SplittingNonlinear}
\begin{equation}
\begin{gathered}
\tfrac{\dd}{\dd t} \, u_1(t) = \alpha F_1\big(u_1(t)\big)\,, \quad \nE_{\tau, \alpha F_1}\big(u_1(t_n)\big) = u_1(t_n + \tau)\,, \\
\tfrac{\dd}{\dd t} \, u_2(t) = \beta F_2\big(u_2(t)\big)\,, \quad \nE_{\tau, \beta F_2}\big(u_2(t_n)\big) = u_2(t_n + \tau)\,, \\
\alpha, \beta \in \CC\,, \quad t \in [t_n, t_n + \tau]\,, 
\end{gathered}
\end{equation}  
so that a higher-order splitting method applied to~\eqref{eq:IVPNonlinear} reads as 
\begin{equation}
\begin{gathered}
u_{n+1} = \big(\nE_{\tau, b_s F_2} \circ \nE_{\tau, a_s F_1} \circ \dots \circ \nE_{\tau, b_1 F_2} \circ \nE_{\tau, a_1 F_1}\big)(u_n) \approx u(t_{n+1})\,, \\
n \in \{0, 1, \dots, N-1\}\,.
\end{gathered}
\end{equation}  
\end{subequations}

\MyParagraph{Schemes}
In view of numerical comparisons, we introduce the coefficients of the first-order Lie--Trotter splitting method 
\begin{subequations}
\label{eq:Coefficients}
\begin{equation}
\label{eq:CoefficientsLie}
s = 1\,, \quad a_1 = 1\,, \quad b_1 = 1\,, 
\end{equation}
and the second-order Strang splitting method 
\begin{equation}
\label{eq:CoefficientsStrang}
s = 2\,, \quad a_1 = 0\,, \quad a_2 = 1\,, \quad b_1 = \tfrac{1}{2} = b_2\,.
\end{equation}
{The well-known fourth-order splitting method by \textsc{Yoshida}~\cite{Yoshida1990} involves four stages
\begin{equation}
\label{eq:CoefficientsYoshida}
\begin{gathered}
s = 4\,, \quad a_1 = 0\,, \quad a_2 = 1 - 2 \, b_2 = a_4\,, \quad a_3 = 4 \, b_2 - 1\,, \\
b_1 = \tfrac{1}{2} - b_2 = b_4\,, \quad b_2 = \tfrac{1}{6} \, \big(1 - \sqrt[3]{2} - \tfrac{1}{2} \sqrt[3]{4}\big) = b_3\,.
\end{gathered}
\end{equation}
Reconsidering its construction based on a triple jump composition of the Strang splitting method, a corresponding fourth-order splitting method with complex coefficients is obtained 
\begin{equation}
\label{eq:CoefficientsYoshidaComplex}
\begin{gathered}
s = 4\,, \quad a_1 = 0\,, \quad a_2 = 1 - 2 \, b_2 = a_4\,, \quad a_3 = 4 \, b_2 - 1\,, \\
b_1 = \tfrac{1}{2} - b_2 = b_4\,, \quad b_2 = \tfrac{1}{6} \, \big(1 + \tfrac{1}{2} \sqrt[3]{2} + \tfrac{1}{4} \sqrt[3]{4}\big) + \ii \, \tfrac{\sqrt{3}}{12} \, \big(\tfrac{1}{2} \sqrt[3]{4} - \sqrt[3]{2}\big) = b_3\,, 
\end{gathered}
\end{equation} 
see also~\cite{BlanesCasasEscorihuela2022,BlanesCasasGonzalezThalhammer2023b,HairerLubichWanner2006}.
}
\end{subequations}

\MyParagraph{Stability issues}
{In connection with the time integration of dissipative systems and parabolic equations as well as the imaginary time propagation of Schr{\"o}dinger equations, operator splitting methods are subject to additional stability constraints. 
In order to explain this matter of fact, we recall that the application of a splitting method with real coefficients $(a_j, b_j)_{j=1}^s$ to the Gross--Pitaevskii equation~\eqref{eq:GPE} in imaginary time and to the parabolic equation~\eqref{eq:Parabolic} involves the subproblems} 
\begin{equation*}
\partial_t U(x, t) = a_j \, \Delta U(x, t)\,, \quad (x, t) \in \Omega \times [t_n, t_n + \tau]\,, \quad j \in \{1, \dots, s\}\,.
\end{equation*}
Evidently, requiring well-posedness of these subproblems or stability of the resulting splitting method, respectively, implies 
\begin{equation*}
a_j \geq 0\,, \quad j \in \{1, \dots, s\}\,.
\end{equation*}
This positivity condition, however, excludes higher-order schemes, since any splitting method that exceeds a second-order barrier necessarily comprises negative coefficients, see for example~\cite{BlanesCasas2005,Sheng1989,Suzuki1991}.
Specifically, this holds true for the fourth-order splitting method by \textsc{Yoshida} 
\begin{equation*}
a_3 \approx - \, 1.7 < 0\,, 
\end{equation*}
see~\eqref{eq:CoefficientsYoshida}.
A feasible remedy to this issue is the design of splitting methods with complex coefficients $(a_j, b_j)_{j=1}^s$ such that 
\begin{equation*}
\Re(a_j)  \geq 0\,, \quad j \in \{1, \dots, s\}\,, 
\end{equation*}
see, e.g., \cite{CastellaChartierDecombesVilmart2009,HansenOstermann2009}.
The fourth-order scheme~\eqref{eq:CoefficientsYoshidaComplex} indeed fulfills these constraints
\begin{equation*}
a_1 = 0\,, \quad \Re(a_2) = \Re(a_4) \approx 0.3 > 0\,, \quad \Re(a_3) \approx 0.4 > 0\,.
\end{equation*}
For further considerations in the context of the imaginary time propagation of the linear Schr{\"o}dinger equation by complex splitting methods, we refer to~\cite{BaderBlanesCasas2013}.

{\MyParagraph{Modified potential splitting methods}}
Reviving former work by \textsc{Ruth} and \textsc{Suzuki}, see for instance~\cite{Ruth1983,Suzuki1991,Suzuki1995}, a favourable alternative to standard operator splitting methods was proposed by \textsc{Chin}.
In a seminal work~\cite{Chin1997}, he developed a famous fourth-order scheme of splitting type that comprises positive coefficients and hence overcomes the second-order barrier for standard splitting methods with real coefficients. 
Expressed in his own words, the basic idea is to incorporate \emph{an additional higher order composite operator} so that the implementation of \emph{one algorithm requires only one evaluation of the force and one evaluation of the force and its gradient}.
For linear evolution equations~\eqref{eq:IVPLinear}, the resulting scheme takes the form
\begin{equation}
\label{eq:ModifiedSplittingLinear}
\begin{gathered}
u_{n+1} = \ee^{\frac{1}{6} \tau B} \, \ee^{\frac{1}{2} \tau A} \, \ee^{\frac{2}{3} \tau B - \frac{1}{72} \tau^3 [B,[B,A]]} \, \ee^{\frac{1}{2} \tau A} \, \ee^{\frac{1}{6} \tau B} \, u_n \approx u(t_{n+1})\,, \\
n \in \{0, 1, \dots, N-1\}\,.
\end{gathered}
\end{equation}
Here, the iterated commutator of complex matrices $A, B \in \CC^{M \times M}$ is given by 
\begin{subequations}
\begin{equation}
\label{eq:IteratedCommutatorLinear}
[B, A] = B A - A B\,, \quad \big[B,[B,A]\big] = B^2 A - 2 \, B A B + A B^2\,,
\end{equation}
see also~\eqref{eq:BBAIntroduction}.
More generally, for linear differential and multiplication operators 
\begin{equation}
\label{eq:B}
\begin{gathered}
(A \, v)(x) = c \, \Delta \, v(x)\,, \\
(B \, v)(x) = \bar{c} \, V(x) \, v(x)\,, \\
x \in \Omega\,, \quad c \in \CC\,,
\end{gathered}
\end{equation}
retained from~\eqref{eq:GPE} and~\eqref{eq:Parabolic} for the special choice $\vartheta = 0$, 
a straightforward calculation yields a linear multiplication operator that depends on the Euclidean norm of the gradient of the potential 
\begin{equation}
\label{eq:BBA}
\begin{split}
&\big[B,[B,A]\big] v(x) \\
&\quad= \bar{c} \, \abs{c}^2 \, \bigg(\big(V(x)\big)^2 \Delta \, v(x) - 2 \, V(x) \, \Delta \big(V(x) \, v(x)\big) \\
&\quad\qquad + \Delta \Big(\big(V(x)\big)^2 \, v(x)\Big)\bigg) \\
&\quad= 2 \, \bar{c} \, \abs{c}^2 \, \big(\nabla V(x)\big)^T \, \nabla V(x) \, v(x)\,, \quad x \in \Omega\,.
\end{split}
\end{equation}
\end{subequations}
This explains the common notion \emph{force-gradient operator splitting method} or \emph{modified potential operator splitting method} for the scheme~\eqref{eq:ModifiedSplittingLinear} and related splitting methods in the context of classical or quantum many-body problems and beyond.
More recent contributions that exploit~\eqref{eq:ModifiedSplittingLinear} for linear ordinary and partial differential equations are, e.g., \cite{AuerKrotscheckChin2001,Chin2005,OmelyanMryglodFolk2002,OmelyanMryglodFolk2003}.

{\MyParagraph{Generalisations to specific classes of nonlinear evolution equations}
We point out that the operator~$B$ defined in~\eqref{eq:B} and the iterated commutator~\eqref{eq:BBA} are of the same nature.
This structural similarity explains the efficiency of \textsc{Chin}’s scheme~\eqref{eq:ModifiedSplittingLinear} for partial differential equations comprising the Laplacian and a potential acting as multiplication operator.
Anticipating the detailed expositions in Sections~\ref{sec:Section3} and~\ref{sec:Section4}, we stress that the generalisations to the time-dependent Gross--Pitaevskii equation~\eqref{eq:GPE} and its parabolic analogue~\eqref{eq:Parabolic} maintain this feature, even though the iterated commutators are more involved.}
\section{Modified operator splitting method}
\label{sec:Section3}
\MyParagraph{Formal generalisation}
{Our educated guess to formally generalise the modified potential operator splitting method~\eqref{eq:ModifiedSplittingLinear} to the significantly more involved case of a nonlinear evolution equation~\eqref{eq:IVPNonlinear} is
\begin{subequations}
\label{eq:ModifiedSplittingNonlinear}
\begin{equation}
\begin{gathered}
u_{n+1}
= \big(\nE_{\tau, \frac{1}{6} F_2} \circ \nE_{\tau, \frac{1}{2} F_1} \circ \nE_{\tau, \frac{2}{3} F_2 - \frac{1}{72} \tau^2 G_2} \circ \nE_{\tau, \frac{1}{2} F_1} \circ \nE_{\tau, \frac{1}{6} F_2}\big) \, u_n \approx u(t_{n+1})\,, \\
n \in \{0, 1, \dots, N-1\}\,. 
\end{gathered}
\end{equation}
That is, we replace the matrices~$A$ and~$B$ by the nonlinear operators~$F_1$ and~$F_2$ as well as the iterated commutator of matrices by the following operator 
\begin{equation}
\label{eq:IteratedCommutator}
\begin{gathered}
G_1(v) = F_2'(v) \, F_1(v) - F_1'(v) \, F_2(v)\,, \\
G_2(v) = F_2'(v) \, G_1(v) - G_1'(v) \, F_2(v)\,, 
\end{gathered}
\end{equation}  
\end{subequations}
see~\eqref{eq:G2Introduction} and~\eqref{eq:BBAIntroduction}. 
We determine the G{\^a}teaux derivatives generalising directional derivatives through 
\begin{equation}
\label{eq:Gateaux}
H'(v) \, w = \lim_{\varepsilon \to 0} \tfrac{1}{\varepsilon} \, \big(H(v + \varepsilon \, w) - H(v)\big)\,, 
\end{equation}
see also \cite{HillePhillips1974}, and presume well-definedness of the arising operators on suitably chosen domains.}
According to~\eqref{eq:SplittingNonlinear}, the decisive operator is associated with a nonlinear evolution equation that comprises the time increment as parameter 
\begin{equation}
\label{eq:EquationF2G2}  
\begin{gathered}
\tfrac{\dd}{\dd t} \, \widetilde{u}_2(t) = \beta_1 F_2\big(\widetilde{u}_2(t)\big) + \beta_2 \, \tau^2 \, G_2\big(\widetilde{u}_2(t)\big)\,, \quad t \in [t_n, t_n + \tau]\,, \\
\nE_{\tau, \, \beta_1 F_2 + \, \beta_2 \tau^2 G_2}\big(\widetilde{u}_2(t_n)\big) = \widetilde{u}_2(t_n + \tau)\,, \quad \beta_1, \beta_2 \in \RR\,.
\end{gathered}
\end{equation}  

\MyParagraph{Specification and implementation}
In the subsequent sections, we complete the remaining tasks. 
We first specify the iterated commutator~\eqref{eq:IteratedCommutator} for the time-dependent Gross--Pitaevskii equation~\eqref{eq:GPE} and contrast it to the result obtained for the parabolic counterpart~\eqref{eq:Parabolic}.
Then, we deduce an invariance principle that has a substantial impact on the efficiency of the modified operator splitting method~\eqref{eq:ModifiedSplittingNonlinear} when applied to the Gross--Pitaevskii equation~\eqref{eq:GPE}.
Implementation issues as well as strategies to reduce the computational cost for parabolic equations are finally discussed in Section~\ref{sec:Section6}.
\section{Iterated commutators}
\label{sec:Section4}
Generally speaking, the appropriate framework for the extension of iterated commutators for matrices or linear operators, see~\eqref{eq:IteratedCommutatorLinear}, to nonlinear operators is provided by the formal calculus of Lie-derivatives.
In the present work, we concretise and verify the heuristic characterisation~\eqref{eq:ModifiedSplittingNonlinear} for relevant applications, the time-dependent Gross--Pitaevskii equation~\eqref{eq:GPE} and the related parabolic equation~\eqref{eq:Parabolic}. 
In the context of Schr{\"o}dinger equations, the arising functions $v, w: \Omega \to \CC$ take complex values, whereas it suffices to consider real-valued functions $v, w: \Omega \to \RR$ for parabolic problems.
In both cases, suitable regularity requirements apply.
For notational simplicity, we omit the dependence on the spatial variable. 

\MyParagraph{Derivatives}
The G{\^a}teaux derivatives of the linear differential and nonlinear multiplication operators defined in~\eqref{eq:GPE} and~\eqref{eq:Parabolic} are given by 
\begin{equation*}  
\begin{gathered}
F_1(v) = c \, \Delta v\,, \\
F_1'(v) \, w = c \, \Delta w\,, \\
F_2(v) = \bar{c} \, \big(V + \vartheta \, \abs{v}^2\big) \, v = \bar{c} \, \big(V v + \vartheta \, v^2 \, \bar{v}\big)\,, \\
F_2'(v) \, w = \bar{c} \, \big(V w + 2 \, \vartheta \, \bar{v} \, v \, w + \vartheta \, v^2 \, \overline{w}\big)\,,
\end{gathered}
\end{equation*}  
see also~\eqref{eq:Gateaux}.

\MyParagraph{First commutators}
On the one hand, by performing differentiation twice, we obtain 
\begin{equation*}
\begin{split}
&F_1'(v) \, F_2(v) \\
&\quad = \abs{c}^2 \, \Delta \big(V v + \vartheta \, v^2 \, \bar{v}\big) \\
&\quad = \abs{c}^2 \, \big(\Delta V v 
+ 2 \, (\nabla V)^T \, \nabla v
+ V \Delta v\big) \\
&\quad\quad\, + \abs{c}^2 \, \vartheta \, \big(2 \, \Delta v \, \bar{v} \, v 
+ \Delta \bar{v} \, v^2
+ 2 \, (\nabla v)^T \, \nabla v \, \bar{v}
+ 4 \, (\nabla v)^T \, \nabla \bar{v} \, v\big)\,.
\end{split}
\end{equation*}
On the other hand, a simple replacement yields 
\begin{equation*}
\begin{split}
&F_2'(v) \, F_1(v) = \abs{c}^2 \, \big(V \Delta v + 2 \, \vartheta \, \bar{v} \, v \, \Delta v\big) + \bar{c}^2 \, \vartheta \, v^2 \Delta \bar{v}\,. 
\end{split}
\end{equation*}
As a consequence, the difference is given by 
\begin{equation*}
\begin{split}
&G_1(v) = F_2'(v) \, F_1(v) - F_1'(v) \, F_2(v) \\  
&\quad = - \, \abs{c}^2 \, \big(\Delta V v + 2 \, (\nabla V)^T \, \nabla v\big) \\
&\quad\quad\, + \big(\bar{c}^2 - \abs{c}^2\big) \, \vartheta \, \Delta \bar{v} \, v^2
- 2 \, \abs{c}^2 \, \vartheta \, \big((\nabla v)^T \, \nabla v \, \bar{v}
+ 2 \, (\nabla v)^T \, \nabla \bar{v} \, v\big)\,.
\end{split}
\end{equation*}
For the parabolic case~\eqref{eq:Parabolic} with $c = 1$ this implies
\begin{equation*}
G_1(v) = - \, \Delta V v - 2 \, (\nabla V)^T \, \nabla v - 6 \, \vartheta \, (\nabla v)^T \, \nabla v \, v\,, 
\end{equation*}
and the analogous result for the Schr{\"o}dinger case~\eqref{eq:GPE} with $c = \ii$ is 
\begin{equation*}
G_1(v) = - \, \Delta V v - 2 \, (\nabla V)^T \, \nabla v - 2 \, \vartheta \, \big(\Delta \bar{v} \, v^2
+ (\nabla v)^T \, \nabla v \, \bar{v}
+ 2 \, (\nabla v)^T \, \nabla \bar{v} \, v\big)\,.
\end{equation*}

\MyParagraph{Iterated commutators}
The iterated commutator associated with the parabolic equation~\eqref{eq:Parabolic} results from straighforward but lengthy calculations
\begin{equation}
\label{eq:G2Parabolic}
\begin{gathered}
G_2(v) = 2 \, \big((\nabla V)^T (\nabla V) + \vartheta \, \widetilde{G}_2(v)\big) \, v\,, \\
\widetilde{G}_2(v) = - \, \Delta V \, v^2
+ 6 \, (\nabla V)^T (\nabla v) \, v
+ 6 \, \big(V + 2 \, \vartheta \, v^2\big) \, (\nabla v)^T (\nabla v)\,, 
\end{gathered}
\end{equation}
and for the Gross--Pitaevskii equation~\eqref{eq:GPE}, we instead arrive at 
\begin{equation}
\label{eq:G2GPE}
\begin{gathered}
G_2(v) = - \, 2 \, \ii \, \Big((\nabla V)^T (\nabla V) - 2 \, \vartheta \, \big(\widetilde{G}_{21}(v) + \vartheta \, \widetilde{G}_{22}(v)\big)\Big) \, v\,, \\  
\widetilde{G}_{21}(v) = \abs{v}^2 \, \Delta V\,, \\
\widetilde{G}_{22}(v) = \abs{v}^2 \, \big(2 \, \Re(\bar{v} \, \Delta v) + 3 \, (\nabla \bar{v})^T (\nabla v)\big) + \Re\big(\bar{v}^2 \, (\nabla v)^T (\nabla v)\big)\,.
\end{gathered}
\end{equation}
For the special case of linear multiplication operators, i.e. $\vartheta = 0$, we indeed recover~\eqref{eq:BBA}.
It is also noteworthy that the operator in~\eqref{eq:G2GPE} involves~$\Delta v$ and hence implies stronger regularity requirements on~$v$ compared to the operator in~\eqref{eq:G2Parabolic}, which only comprises the gradient~$\nabla v$. 
\section{Invariance principle}
\label{sec:Section5}
{\MyParagraph{Invariance principle}
For the purpose of illustration, we introduce the nonlinear Schr{\"o}dinger equation 
\begin{equation*}
\begin{cases}
\ii \, \partial_t \Psi(x, t) = - \, \Delta \Psi(x, t) + \abs{\Psi(x, t)}^2 \, \Psi(x, t)\,, \\
\Psi(x, 0) = \Psi_0(x)\,, \quad (x, t) \in \Omega \times [0, T]\,, 
\end{cases}
\end{equation*}  
retained from~\eqref{eq:GPE} for vanishing potential~$V = 0$ and normalised constant $\vartheta = 1$.   
Any standard operator splitting method~\eqref{eq:SplittingNonlinear} relies on the time integration of the subproblem involving the Laplacian 
\begin{equation*}
\begin{cases}
\ii \, \partial_t \Psi_1(x, t) = - \, \alpha \, \Delta \Psi_1(x, t)\,, \quad \alpha \in \RR\,, \\
\Psi_1(x, t_n) \text{ given}\,, \quad (x, t) \in \Omega \times [t_n, t_n + \tau]\,, 
\end{cases}
\end{equation*}
which we perform by means of a Fourier spectral space discretisation, see Section~\ref{sec:Section6}.
Moreover, we make use of the fact that the subproblem comprising the cubic nonlinearity 
\begin{subequations}
\label{eq:NLS2}
\begin{equation}
\begin{cases}
\ii \, \partial_t \Psi_2(x, t) = \beta_1 \, \abs{\Psi_2(x, t)}^2 \, \Psi_2(x, t)\,, \quad \beta_1 \in \RR\,, \\
\Psi_2(x, t_n) \text{ given}\,, \quad (x, t) \in \Omega \times [t_n, t_n + \tau]\,, 
\end{cases}
\end{equation}
satisfies the fundamental invariance principle 
\begin{equation*}
\abs{\Psi_2(x, t)}^2 = \abs{\Psi_2(x, t_n)}^2\,, \quad (x, t) \in \Omega \times [t_n, t_n + \tau]\,.
\end{equation*}
This identity follows by means of differentiation with respect to time and substitution of the time derivative by the nonlinearity 
\begin{equation*}
\begin{gathered}  
\partial_t \, \abs{\Psi_2(x, t)}^2
= 2 \, \Re\big(\overline{\Psi_2(x, t)} \, \partial_t \Psi_2(x, t)\big)
= - \, 2 \, \beta_1 \, \Re\big(\ii \, \abs{\Psi_2(x, t)}^4\big)
= 0\,, \\
(x, t) \in \Omega \times [t_n, t_n + \tau]\,.
\end{gathered}
\end{equation*}
It implies that the exact solution to the nonlinear subproblem satisfies the reduced problem
\begin{equation}
\begin{cases}
\ii \, \partial_t \Psi_2(x, t) = \beta_1 \, \abs{\Psi_2(x, t_n)}^2 \, \Psi_2(x, t)\,, \quad \beta_1 \in \RR\,, \\
\Psi_2(x, t_n) \text{ given}\,, \quad (x, t) \in \Omega \times [t_n, t_n + \tau]\,, 
\end{cases}
\end{equation}
and hence can be determined by pointwise multiplication 
\begin{equation}
\Psi_2(x, t) = \ee^{- \, \ii \, \beta_1 \, \abs{\Psi_2(x, t_n)}^2} \, \Psi_2(x, t_n)\,, \quad (x, t) \in \Omega \times [t_n, t_n + \tau]\,.
\end{equation}
\end{subequations}
For the numerical realisation, equidistant grid points are chosen accordingly to the Fourier spectral space discretisation, see also Section~\ref{sec:Section6}.

\MyParagraph{Generalisation}
In this section, we establish a substantial extension of the above stated result~\eqref{eq:NLS2} that has important implications concerning the efficient implementation of the modified operator splitting method~\eqref{eq:ModifiedSplittingNonlinear} for the time-dependent Gross--Pitaevskii equation~\eqref{eq:GPE}, see also~\eqref{eq:EquationF2G2} and~\eqref{eq:G2GPE}.  
Furthermore, the obtained invariance principle for the evolution operator
\begin{equation*}
\nE_{\tau, \, \beta_1 F_2 + \, \beta_2 \tau^2 G_2}\,, \quad \tau, \beta_1, \beta_2 \in \RR\,, 
\end{equation*}
is connected with the significance of the modified operator splitting method as a geometric integrator, see~\cite{HairerLubichWanner2006,Iserles2008} and references given therein.
}

\MyParagraph{Notation}
With regard to a compact formulation as abstract evolution equation, we again omit the dependence of the potential and a regular complex-valued function $v: \Omega \to \CC$ on the spatial variable.
Besides, for accomplishing relations of the form 
\begin{equation*}
F_2(v) = \bar{c} \, f_1(v) \, v\,, \quad G_2(v) = \bar{c} \, f_2(v) \, v\,, 
\end{equation*}
it is convenient to introduce the abbreviations 
\begin{equation}
\label{eq:fDefinition}  
\begin{gathered}
f(v) = \beta_1 \, f_1(v) + \beta_2 \, \tau^2 f_2(v)\,, \\
f_1(v) = V + \vartheta \, g_1(v)\,, \quad f_2(v) = 2 \, (\nabla V)^T (\nabla V) - 4 \, \vartheta \, g_6(v)\,, \\
g_1(v) = \abs{v}^2\,, \quad 
g_2(v) = \Re(\bar{v} \, \Delta v)\,, \\
g_3(v) = (\nabla \bar{v})^T (\nabla v)\,, \quad
g_4(v) = \Re\big(\bar{v}^2 \, (\nabla v)^T (\nabla v)\big)\,, \\
g_5(v) = \vartheta \, \big(2 \, g_2(v) + 3 \, g_3(v)\big)\,, \quad
g_6(v) = g_1(v) \, \big(\Delta V + g_5(v)\big) + \vartheta \, g_4(v)\,.
\end{gathered}
\end{equation}

\MyParagraph{Theorem (Invariance principle)}
The solution to the subproblem 
\begin{equation*}
\begin{cases}
\tfrac{\dd}{\dd t} \, \psi(t) = - \, \ii \, f\big(\psi(t)\big) \, \psi(t)\,, \\
\psi(0) = \psi_0\,, \quad t \in [0, \tau]\,,  
\end{cases}
\end{equation*}
with defining function introduced in~\eqref{eq:fDefinition} satisfies the invariance principle 
\begin{equation*}
f\big(\psi(t)\big) = f(\psi_0)\,, \quad t \in [0, \tau]\,.
\end{equation*}  
\MyParagraph{Proof}
In order to demonstrate that the invariance principle holds, we determine the G{\^a}teaux derivatives of the defining functions 
\begin{equation*}
\begin{gathered}
f'(v) \, w = \beta_1 \, f_1'(v) \, w + \beta_2 \, \tau^2 f_2'(v) \, w\,, \\
f_1'(v) \, w = \vartheta \, g_1'(v) \, w\,, \quad
f_2'(v) \, w = - \, 4 \, \vartheta \, g_6'(v) \, w\,, \\
g_1'(v) \, w = 2 \, \Re(\bar{v} \, w)\,, \quad
g_2'(v) \, w = \Re(\Delta \bar{v} \, w) + \Re(\bar{v} \, \Delta w)\,, \\
g_3'(v) \, w = 2 \, \Re\big((\nabla \bar{v})^T (\nabla w)\big)\,, \\
g_4'(v) \, w = 2 \, \Re\big(\bar{v} \, \overline{w} \, (\nabla v)^T (\nabla v)\big) + 2 \, \Re\big(\bar{v}^2 \, (\nabla v)^T (\nabla w)\big)\,, \\
g_5'(v) \, w = 2 \, \vartheta \, g_2'(v) \, w + 3 \, \vartheta \, g_3'(v) \, w\,, \\
g_6'(v) \, w = \big(\Delta V + g_5(v)\big) \, g_1'(v) \, w + g_1(v) \, g_5'(v) \, w + \vartheta \, g_4'(v) \, w\,, \\
\end{gathered}
\end{equation*}
where suitable regularity requirements apply to $v, w: \Omega \to \CC$.
Observing that the potential~$V$ and the basic components $g_1, g_2, g_3, g_4$ define real-valued functions, we have  
\begin{equation*}
\begin{gathered}
\big(g_j(v)\big)(x) \in \RR\,, \quad j \in \{1, \dots, 6\}\,, \\
\big(f_k(v)\big)(x) \in \RR\,, \quad k \in \{1,2\}\,, \quad \big(f(v)\big)(x) \in \RR\,, \quad x \in \Omega\,.
\end{gathered}
\end{equation*}
Evidently, this implies that the following composition vanishes 
\begin{equation*}
g_1'(v) \, \big(\ii \, f(v) \, v\big) = 2 \, \Re\big(\ii \, \abs{v}^2 \, f(v)\big) = 0\,.
\end{equation*}
Certain contributions originating from the iterated commutator, however, require a closer examination, namely
\begin{equation*}
\begin{split}
&g_2'(v) \, \big(\ii \, f(v) \, v\big) \\
&\quad = \Re\Big(\ii \, f(v) \, \Re\big(\bar{v} \, \Delta v\big)\Big)
+ \Re\big(\ii \, \abs{v}^2 \Delta f(v)\big) 
+ 2 \, \Re\Big(\ii \, \bar{v} \, \big(\nabla f(v)\big)^T \, \nabla v\Big) \\
&\quad = 2 \, \Re\Big(\ii \, \bar{v} \, \big(\nabla f(v)\big)^T \nabla v\Big)\,, \\
&g_3'(v) \, \big(\ii \, f(v) \, v\big) \\
&\quad = 2 \, \Re\Big(\ii \, v \, \big(\nabla f(v)\big)^T \, \nabla \bar{v}\Big) + 2 \, \Re\Big(\ii \, f(v) \, (\nabla \bar{v})^T \, \nabla v\Big) \\
&\quad = 2 \, \Re\Big(\ii \, v \, \big(\nabla f(v)\big)^T \, \nabla \bar{v}\Big)\,, \\
&g_4'(v) \, \big(\ii \, f(v) \, v\big) \\
&\quad = - \, 2 \, \Re\big(\ii \, f(v) \, \bar{v}^2 \, (\nabla v)^T \, \nabla v\big)
+ 2 \, \Re\Big(\ii \, \abs{v}^2 \, \bar{v} \, \big(\nabla f(v)\big)^T \nabla v\Big) \\
&\qquad + 2 \, \Re\big(\ii \, f(v) \, \bar{v}^2 \, (\nabla v)^T \, \nabla v\big) \\
&\quad = 2 \, \Re\Big(\ii \, \abs{v}^2 \, \bar{v} \, \big(\nabla f(v)\big)^T \nabla v\Big)\,.
\end{split}
\end{equation*}
On the basis of these identities, we conclude 
\begin{equation*}
\begin{split}
&\Big(g_1(v) \, \big(2 \, g_2'(v) + 3 \, g_3'(v)\big) + g_4'(v)\big)\Big) \, \big(\ii \, f(v) \, v\big) \\
&\quad = 12 \, \Re\Big(\ii \, \abs{v}^2 \, \big(\nabla f(v)\big)^T \Re\big(\bar{v} \, \nabla v\big)\Big) = 0\,.
\end{split}
\end{equation*}
This proves that any composition of the special form 
\begin{equation*}
\begin{split}
&f'(v) \, \big(\ii \, f(v) \, v\big) \\
&\quad = \vartheta \, \Big(\beta_1 - 4 \, \beta_2 \, \tau^2 \, \big(\Delta V + g_5(v)\big)\Big) \, g_1'(v) \, \big(\ii \, f(v) \, v\big) \\
&\qquad\, - 4 \, \vartheta^2 \beta_2 \, \tau^2 \Big(g_1(v) \, \big(2 \, g_2'(v) + 3 \, g_3'(v)\big) + g_4'(v)\Big) \, \big(\ii \, f(v) \, v\big) \\
&\quad = 0
\end{split}
\end{equation*}
vanishes.
As a consequence, the time-derivative of the decisive function is equal to zero 
\begin{equation*}
\tfrac{\dd}{\dd t} \, f\big(\psi(t)\big) = f'\big(\psi(t)\big) \, \tfrac{\dd}{\dd t} \, \psi(t) = - \, f'\big(\psi(t)\big) \, \Big(\ii \, f\big(\psi(t)\big) \, \psi(t)\Big) = 0\,, \quad t \in [0, \tau]\,, 
\end{equation*}
and hence, the desired identity follows 
\begin{equation*}
f\big(\psi(t)\big) = f(\psi_0)\,, \quad t \in [0, \tau]\,.
\end{equation*}  

\MyParagraph{Summary}
{In view of Section~\ref{sec:Section6}, we summarise our considerations for the modified operator splitting method~\eqref{eq:ModifiedSplittingNonlinear} applied to the time-dependent Gross--Pitaevskii equation~\eqref{eq:GPE}.
On the one hand, the realisation of
\begin{equation*}
\begin{gathered}
\psi_{n+1}
= \big(\nE_{\tau, \frac{1}{6} F_2} \circ \nE_{\tau, \frac{1}{2} F_1} \circ \nE_{\tau, \frac{2}{3} F_2 - \frac{1}{72} \tau^2 G_2} \circ \nE_{\tau, \frac{1}{2} F_1} \circ \nE_{\tau, \frac{1}{6} F_2}\big) \, \psi_n\,, \\
n \in \{0, 1, \dots, N-1\}\,,
\end{gathered}
\end{equation*}
relies on the numerical integration of the linear Schr{\"o}dinger equation 
\begin{equation*}
\begin{gathered}
\tfrac{\dd}{\dd t} \, \psi(t) = \ii \, \alpha \, \Delta \psi(t)\,, \quad t \in [t_n, t_n + \tau]\,, \quad
\nE_{\tau, \alpha F_1}\big(\psi(t_n)\big) = \psi(t_n + \tau)\,, 
\end{gathered}
\end{equation*}
for $\alpha \in \RR$.
On the other hand, it reduces to the pointwise evaluation of the solution representation
\begin{equation}
\label{eq:EF2G2GPE}  
\nE_{\tau, \, \beta_1 F_2 + \, \beta_2 \tau^2 G_2}(\psi_0)
= \ee^{- \, \ii \, \tau (\beta_1 f_1(\psi_0) + \beta_2 \tau^2 f_2(\psi_0))} \, \psi_0\,, \quad \tau \in \RR\,, 
\end{equation}
for appropriate choices of $\beta_1, \beta_2 \in \RR$, see~\eqref{eq:fDefinition} for the definitions of~$f_1$ and~$f_2$.
We recall that the known identity~\eqref{eq:NLS2} for the particular case $\beta_2 = 0$ is established by a simplified argument.}
\section{Numerical results}
\label{sec:Section6}
In the following, we illustrate the stability and global error behaviour of the novel modified operator splitting method~\eqref{eq:ModifiedSplittingNonlinear} for the time-dependent Gross--Pitaevskii equation~\eqref{eq:GPE} and its parabolic analogue~\eqref{eq:Parabolic}.
The numerical tests, performed in one, two, and three space dimensions, in particular confirm the theoretical considerations of Sections~\ref{sec:Section3} to~\ref{sec:Section5}.
For the purpose of comparison, we in addition include the corresponding results for widely-used standard splitting methods.
Further information on a publicly accessible \textsc{Matlab} code to reproduce Figures~\ref{fig:Figure1}--\ref{fig:Figure4} is found in Appendix~\ref{sec:AppendixMatlabCode}.

\MyParagraph{Implementation}
The practical realisation of standard operator splitting methods such as~\eqref{eq:Coefficients} and of the modified operator splitting method~\eqref{eq:ModifiedSplittingNonlinear}, respectively, requires the time integration of the subproblems involving the Laplacian and the nonlinear multiplication operator. 
In our implementation, we make use of fast Fourier techniques, which are based on the following considerations.
\begin{enumerate}[(i)]
\item
\emph{Space grid.} \; 
With regard to the employed Fourier spectral space discretisation, we replace the underlying unbounded domain by a Cartesian product of sufficiently large intervals 
\begin{equation*}
a = 10\,, \quad x \in [- \, a, a]^d \subset \Omega\,, 
\end{equation*}
and choose the total number of equidistant spatial grid points according to the dimension 
\begin{equation*}
d = 1: \; M = 512\,, \quad d = 2: \; M = 128^2\,, \quad d = 3: \; M = 64^3\,.
\end{equation*}
\item 
\emph{Derivatives.} \;
The iterated commutator arising in the modified operator splitting method~\eqref{eq:ModifiedSplittingNonlinear} for~\eqref{eq:GPE} and~\eqref{eq:Parabolic} involves the gradient~$\nabla V$ and the Laplacian~$\Delta V$ of the space-dependent potential, which we may assume to be known analytically.
Otherwise, we employ the approach described subsequently. 
The numerical computation of the spatial derivatives~$\nabla v$ and~$\Delta v$, where~$v$ represents the current value of the time-discrete solution, is traced back to a fast Fourier transform, pointwise multiplications, and an inverse fast Fourier transform.
Denoting by $(\nF_m)_{m \in \ZZ^d}$ the Fourier functions with periodicity domain $[- \, a, a]^d$, by $(\mu_m)_{m \in \ZZ^d}$ the purely imaginary eigenvalues associated with the first spatial derivatives, and by $(\lambda_m)_{m \in \ZZ^d}$ the corresponding real eigenvalues of the Laplace operator 
\begin{equation*}
\begin{gathered}
\nF_m(x) = (2 \, a)^{- \, \frac{d}{2}} \; \ee^{\, \ii \, \pi \, m_1 \, (\frac{x_1}{a} + 1)} \cdots \ee^{\, \ii \, \pi \, m_d \, (\frac{x_d}{a} + 1)}\,, \\
\nabla \nF_m = \mu_m \, \nF_m\,, \quad \mu_m = \tfrac{\ii \, \pi \, m}{a} \in \CC^{d \times 1}\,, \\ 
\Delta \, \nF_m = \lambda_m \, \nF_m\,, \quad \lambda_m = - \, \tfrac{\pi^2 \abs{m}^2}{a^2} \in \RR\,, \\
m = (m_1, \dots, m_d) \in \ZZ^d\,, \quad x = (x_1, \dots, x_d) \in [- \, a, a]^d\,, 
\end{gathered}  
\end{equation*}  
the following formal representations hold 
\begin{equation*}
\begin{gathered}
v = \sum_{m \in \ZZ^d} v_m \, \nF_m\,, \quad v_m = \int_{[- \, a, a]^d} v(x) \, \nF_{-m}(x) \; \dd x\,, \quad m \in \ZZ^d\,, \\
\nabla v = \sum_{m \in \ZZ^d} \mu_m \, v_m \, \nF_m\,, \quad \Delta v = \sum_{m \in \ZZ^d} \lambda_m \, v_m \, \nF_m\,.
\end{gathered}  
\end{equation*}  
Their realisation relies on a suitable truncation of the infinite index set $\nM \subset \ZZ^d$ so that $\abs{\nM} = M$ and quadrature approximations by the trapezoidal rule. 
\item
\emph{Linear subproblem.} \;
Formally, the solution to the linear subproblem is given by a Fourier series 
\begin{equation*}
\begin{gathered}
\tfrac{\dd}{\dd t} \, u_1(t) = \alpha F_1\big(u(t) \big) = c \, \alpha \, \Delta \, u_1(t) \,, \quad t \in [t_n, t_n + \tau]\,, \quad \alpha \in \RR\,, \\
\nE_{\tau, \alpha F_1}\big(u(t_n)\big) = \nE_{\tau, \alpha F_1}\bigg(\sum_{m \in \ZZ^d} u_{1,m}(t_n) \, \nF_m\bigg) = \sum_{m \in \ZZ^d} \ee^{\, c \, \alpha \, \tau \lambda_m} \, u_{1,m}(t_n) \, \nF_m\,,
\end{gathered}  
\end{equation*}  
see also~\eqref{eq:GPE},~\eqref{eq:Parabolic}, and~\eqref{eq:SplittingNonlinear}.
Again, the application of fast Fourier techniques permits the efficient computation of approximations to the spectral coefficients and the evaluation of finite sums on the equidistant spatial grid points covering the underlying domain
\begin{equation*}
\begin{gathered}
\widetilde{u}_{1,m,n} \approx u_{1,m}(t_n)\,, \quad m \in \nM\,, \\
\sum_{m \in \nM} \ee^{\, c \, \alpha \, \tau \lambda_m} \, \widetilde{u}_{1,m,n} \, \nF_m(x)\,, \quad x \in [- \, a, a]^d \,.
\end{gathered}  
\end{equation*}  
\item
\emph{Nonlinear subproblem (Schr{\"o}dinger equation).} \; 
In the context of the Gross--Pitaevskii equation~\eqref{eq:GPE}, we make use of the fact that the solution to the nonlinear subproblem~\eqref{eq:EquationF2G2} satisfies the invariance principle deduced in Section~\ref{sec:Section5}.
Consequently, it simply remains to evaluate the representation~\eqref{eq:EF2G2GPE} on the equidistant grid.
\item
\emph{Nonlinear subproblem (Parabolic equation).} \; 
In the case of the parabolic equation~\eqref{eq:Parabolic}, we additionally apply an explicit Runge--Kutta method of order four.
Due to the stiffness of the problem, the time stepsize has to be adjusted to the spatial grid width to ensure stability.
Alternative approaches with improved stability properties and reduced computational costs are detailed below. 
\end{enumerate}

\MyParagraph{Numerical results}
In our numerical tests, we perform the time integration of the Gross--Pitaevskii equation~\eqref{eq:GPE} and the parabolic equation~\eqref{eq:Parabolic}, {comparing} the nonlinear case with $\vartheta = 1$ to the simplified linear case with $\vartheta = 0$.
We prescribe the Gaussian-shaped initial state
\begin{equation*}
u_0(x) = \ee^{- \frac{1}{2} \, (x_1^2 + \dots + x_d^2)}\,, \quad x = (x_1, \dots, x_d) \in \RR^d\,,
\end{equation*}
as well as the two polynomial potentials 
\begin{equation*}
\begin{gathered}
V(x) = C_0 \, C_q \sum_{j=1}^{d} x_j^q\,, \\
C_2 = 1\,, \quad C_4 = \tfrac{1}{24}\,, \quad q \in \{2,4\}\,, \quad x = (x_1, \dots, x_d) \in \RR^d\,, 
\end{gathered}
\end{equation*}
for which the {required} first- and second-order derivatives are known analytically.
For the special case of a quadratic potential with prefactor $C_0 \in \{1,- \, 1\}$ chosen accordingly to the type of the equation and $\vartheta = 0$, the knowledge of the exact solution
\begin{equation*}
u(x, t) = \ee^{d \, \bar{c} \, t} \, u_0(x)\,, \quad (x, t) \in \RR^d \times [0, T]\,,
\end{equation*}
permits the validation of the Fourier spectral space discretisation and the time-splitting approach. 
In the general case, we instead determine a numerical reference solution based on a refined time stepsize.
The global errors of the modified operator splitting method~\eqref{eq:ModifiedSplittingNonlinear} at final time $T = 1$, measured in a discrete $L_2$-norm, are compared to those obtained by the standard Lie--Trotter, Strang, and Yoshida splitting methods of non-stiff orders one, two, and four, see also~\eqref{eq:Coefficients}.
The obtained results, displayed in Figures~\ref{fig:Figure1}--\ref{fig:Figure4}, confirm the favourable behaviour of the modified operator splitting method.

{\MyParagraph{Energy conservation}
In order to complement our numerical comparisons regarding the stability and global error behaviour of standard and modified operator splitting methods, we perform the time integration of the one-dimensional Gross--Pitaevskii equation~\eqref{eq:GPE} on the interval~$[0, T]$.
We determine approximations to the values of the energy at equidistant time grid points
\begin{subequations}
\label{eq:Energy}  
\begin{equation}
\begin{split}
E(t_n) &= \int_{\Omega} \big(- \Delta \Psi(x, t_n) + V(x) \, \Psi(x, t_n) + \vartheta \, \abs{\Psi(x, t_n)}^2 \, \Psi(x, t_n)\big) \\
&\qquad\quad \times \overline{\Psi(x, t_n)} \; \dd x\,, \quad n \in \{0, 1, \dots, N\}\,,
\end{split}
\end{equation}
and their deviations with respect to the minimal value
\begin{equation}
E(t_n) - \min\big\{E(t_{\ell}): \ell \in \{0, 1, \dots, N\}\big\}\,, \quad n \in \{0, 1, \dots, N\}\,.
\end{equation}
\end{subequations}
The obtained results confirm the favourable geometric properties of the modified operator splitting method~\eqref{eq:ModifiedSplittingNonlinear}, see Figure~\ref{fig:Figure6}.}

\MyParagraph{Computational cost}
In general, an expedient measure for the computational cost of the modified operator splitting method~\eqref{eq:ModifiedSplittingNonlinear} for the Gross--Pitaevskii equation~\eqref{eq:GPE} and the parabolic equation~\eqref{eq:Parabolic}, respectively, is the number of fast Fourier transforms and their inverses.
Evidently, the numerical solution of two linear subproblems per time step amounts to two fast Fourier transforms and two inverse fast Fourier transforms.
Besides, the time integration of the nonlinear subproblem requires the computation of space derivatives via Fourier transforms. 
\begin{enumerate}[(i)]
\item 
For nonlinear Schr{\"o}dinger equations such as~\eqref{eq:GPE}, the validity of the invariance property permits to significantly reduce the cost related to the evaluation of~$\nE_{\tau, \frac{2}{3} F_2 - \frac{1}{72} \tau^2 G_2}$.
Due to the fact that the spectral coefficients of the current numerical solution~$v$ are available, the computation of gradient~$\nabla v$ and Laplacian~$\Delta v$ results in $d + 1$ inverse fast Fourier transforms.
\item 
For the parabolic equation~\eqref{eq:Parabolic}, a favourable approach is based on the following considerations. 
The presence of the additional factor~$\tau^2$ in connection with the double commutator permits to use an approximation by means of the the second-oder Strang splitting method, that is
\begin{equation*}
\nE_{\tau, \frac{2}{3} F_2 - \frac{1}{72} \tau^2 G_2} \approx \nE_{\frac{1}{2} \tau, \frac{2}{3} F_2} \circ \nE_{\tau, - \frac{1}{72} \tau^2 G_2} \circ \nE_{\frac{1}{2} \tau, \frac{2}{3} F_2}\,.
\end{equation*}
{For linear ordinary differential equations of the form
\begin{equation*}
\tfrac{\dd}{\dd t} \, u(t) = B \, u(t) + \tau^2 \, C \, u(t)\,, \quad t \in [0, \tau]\,, 
\end{equation*}
defined by complex matrices $B, C \in \CC^{M \times M}$, elementary calculations based on series expansions of matrix exponentials confirm that this approach leads to a fourth-order approximation
\begin{equation*}
\ee^{\tau (B + \tau^2 C)} = \ee^{\frac{1}{2} \tau B} \, \ee^{\tau^3 C} \, \ee^{\frac{1}{2} \tau B} + \nO\big(\tau^5\big)\,.
\end{equation*}
The rigorous generalisation of the argument to~\eqref{eq:Parabolic} is part of a future convergence analysis.}
On the one hand, an explicit representation of the evolution operator associated with the nonlinear ordinary differential equation 
\begin{equation*}
\tfrac{\dd}{\dd t} \, u(x, t) = \big(V(x) + \vartheta \, \abs{u(x, t)}^2\big) \, u(x, t)\,, \quad x \in \Omega\,, \quad t \in [0, \tau]\,, 
\end{equation*}
is known.
{More precisely, fixing $x \in \Omega$, the distinction of cases yields 
\begin{equation*}
u(x, t) =
\begin{cases}
\frac{u(x, 0)}{\sqrt{\ee^{- 2 \, t \, V(x)} + \vartheta \, (\ee^{- 2 \, t \, V(x)} - 1)/V(x) \, (u(x, 0))^2}}\,, & V(x) \neq 0\,, \\
\frac{u(x, 0)}{\sqrt{1 - 2 \, t \, \vartheta \, (u(x, 0))^2}}\,, & V(x) = 0\,, \\
\end{cases} \quad t \in [0, \tau]\,.
\end{equation*}  
}%
On the other hand, for the time integration of the nonlinear subproblem involving the double commutator, it suffices to apply the first-order explicit Euler method. 
Numerical tests confirm the enhanced stability, accuracy, and efficiency of the resulting approach, in particular, in case of a fourth-order polynomial potential, see Figure~\ref{fig:Figure5}.
\end{enumerate}
\section{Conclusions}
In the present work, we have introduced a general framework for the extension of \textsc{Chin}’s fourth-order modified potential operator splitting method~\eqref{eq:ModifiedSplittingLinear} to nonlinear evolution equations. 
To the best of our knowledge, this matter is novel and of major interest from theoretical and practical perspectives.

Moreover, we have specified the resulting fourth-order modified operator splitting method~\eqref{eq:ModifiedSplittingNonlinear} for the time-dependent Gross--Pitaevskii equation~\eqref{eq:GPE} and its parabolic counterpart~\eqref{eq:Parabolic}. 
It seems likely that our approach and the drawn conclusions also extend to Schr{\"o}dinger equations involving nonlinear and possibly nonlocal terms of the form 
$g(\abs{\Psi(x, t)}) \, \Psi(x, t)$\,. 

Due to the fact that our numerical tests have confirmed the favourable performance of the proposed fourth-order scheme in comparison with standard real and complex splitting methods of order four, it is natural to extend our considerations in various respects.

Proceeding our recent work~\cite{BlanesCasasGonzalezThalhammer2023} for linear evolution equations of parabolic and Schr{\"o}dinger type, we find it promising to design high-order modified operator splitting methods for nonlinear evolution equations that are optimal with regard to a preselected criterium such as efficiency. 
In contrast to the linear case, we have to take into account the additional costs for the evaluation of second iterated commutators and that higher-order iterated commutators will not vanish, in general, see~\eqref{eq:BBA} as well as~\eqref{eq:G2Parabolic} and~\eqref{eq:G2GPE}. 
Nonetheless, based on the successful strategies for the efficient implementation of the fourth-order scheme~\eqref{eq:ModifiedSplittingNonlinear}, it suggests itself to address an in-depth analysis of similarly structured modified operator splitting methods.

A desirable feature of the proposed fourth-order modified operator splitting method is the positivity of the coefficients.
Concerning the design of high-order schemes, in light of contributions on the linear case, see for instance~\cite{AuzingerHofstaetterKoch2019,Chin2005}, this will necessitate further deliberations on the appropriate format.
Evidently, positivity is intrinsically related to the issue of well-posedness and stability for parabolic problems.
Besides, it affects aliasing effects on truncated space domains and the incorporation of artificial boundary conditions, which is of particular interest in the context of nonlinear Schr{\"o}dinger equations, see for example~\cite{AntoineEtAl2008} and references given therein.

{Furthermore, we intend to carry out a rigorous convergence analysis of modified operator splitting methods applied to the time-dependent Gross--Pitaevskii equation and its parabolic analogue.} 
We point out that the study for the linear case~\cite{Kieri2015} together with the formal calculus of Lie-derivatives provides a guiding principle. 
However, the accomplishment for specific nonlinear evolution equations involving unbounded operators requires careful calculations and investigations, see also~\cite{Thalhammer2008}.
\section*{Acknowledgements}
The authors dedicate this work to \textsc{Arieh Iserles} due to his seminal contributions in the area of numerical analysis and geometric numerical integration.
{The authors are grateful to the two anonymous reviewers for their careful reading of the manuscript and valuable comments.}
This work has been supported by Ministerio de Ciencia e Innovaci{\'o}n (Spain) through projects PID2019-104927GB-C21 and PID2019-104927GB-C22, MCIN/AEI/10.13039/501100011033, ERDF (\emph{A way of making Europe}).
Sergio Blanes and Fernando Casas acknowledge the support of the Conselleria d'Innovaci{\'o}, Universitats, Ci{\`e}ncia i Societat Digital from the Generalitat Valenciana (Spain) through project CIAICO/2021/180.

\begin{figure}[t!]
\begin{center}
\includegraphics[width=6cm]{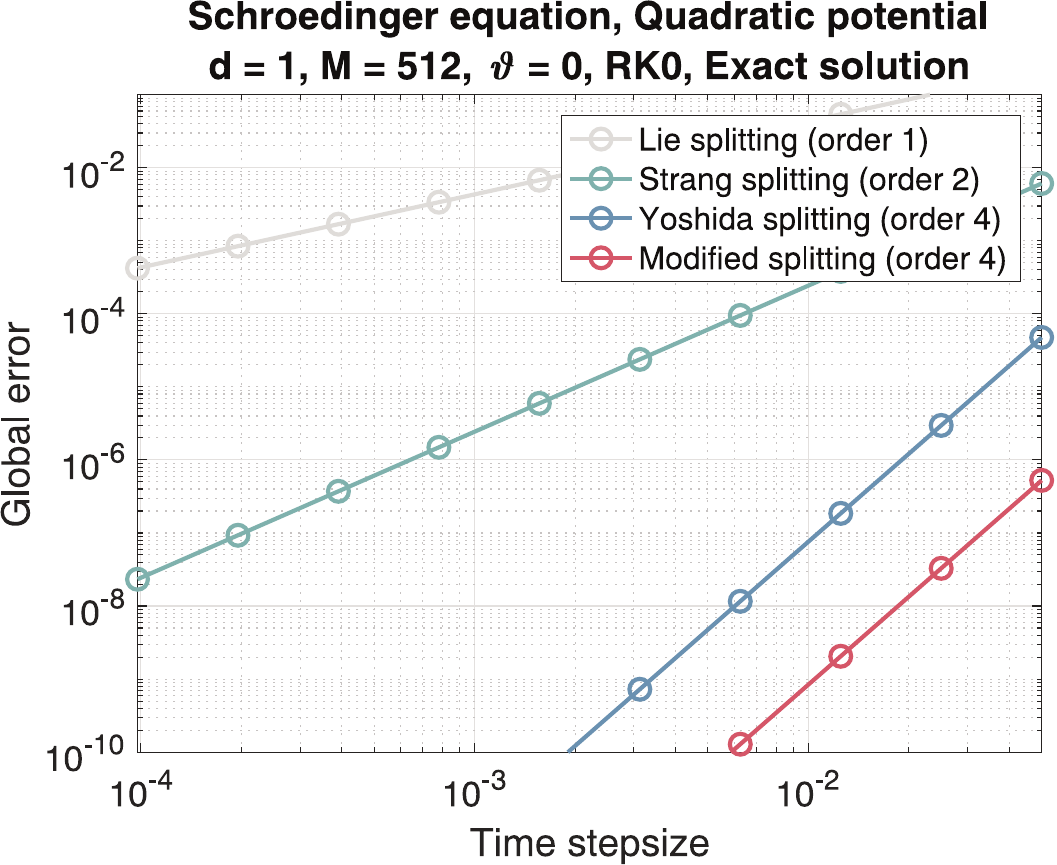} \quad
\includegraphics[width=6cm]{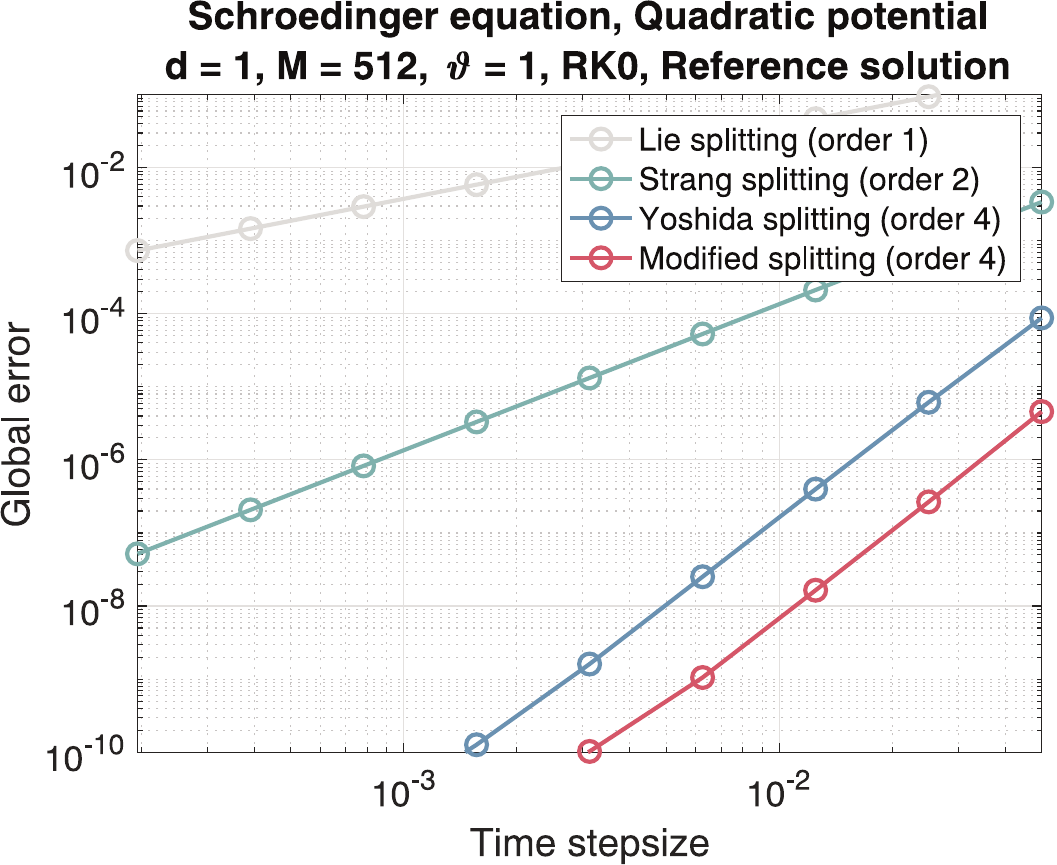} \\[4mm]
\includegraphics[width=6cm]{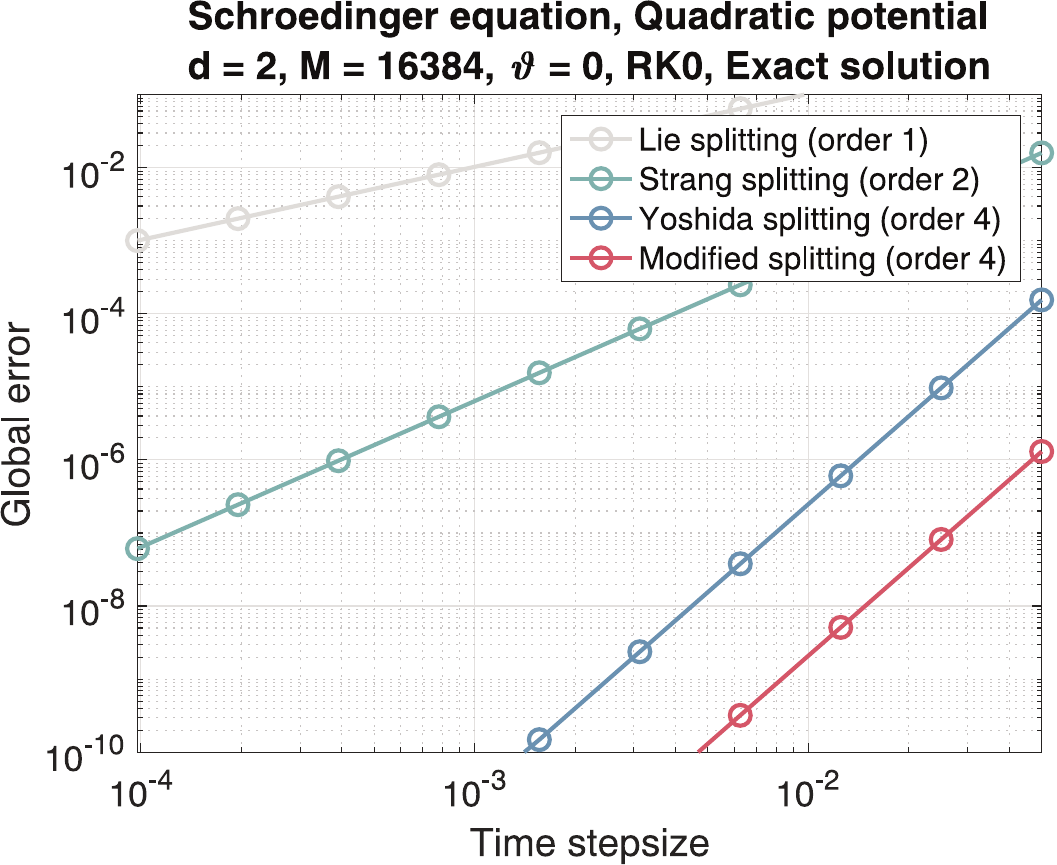} \quad
\includegraphics[width=6cm]{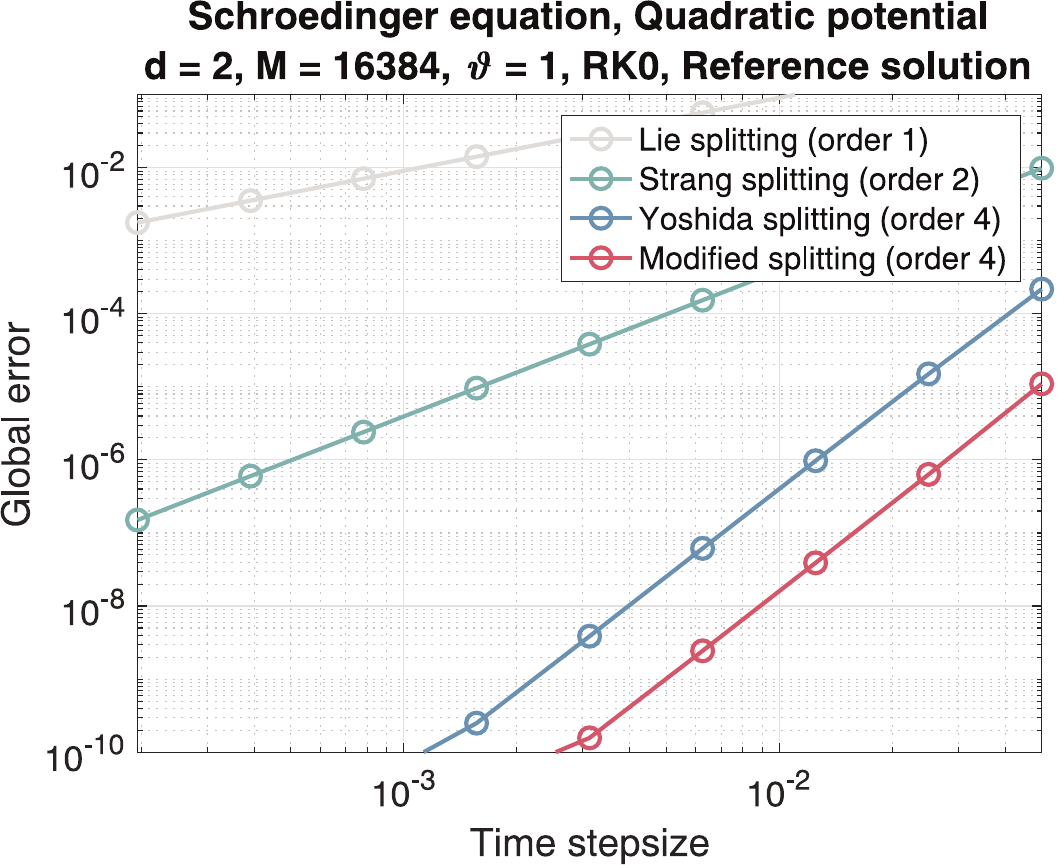} \\[4mm]
\includegraphics[width=6cm]{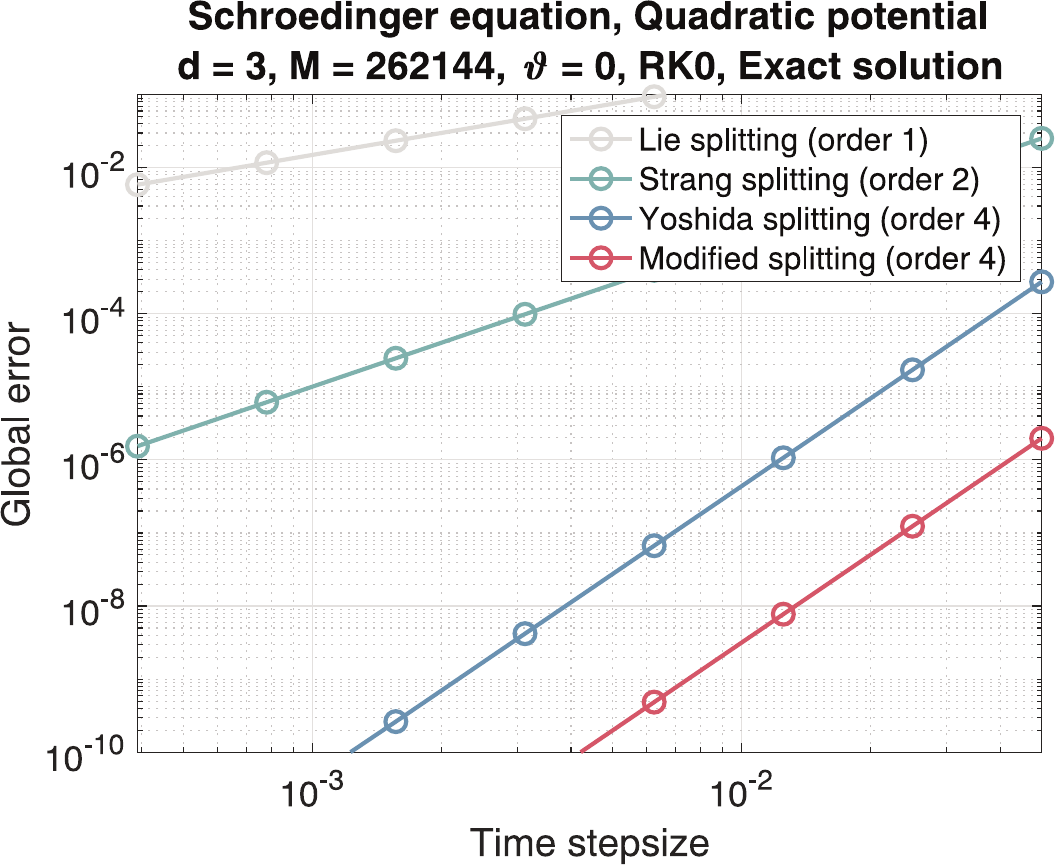} \quad
\includegraphics[width=6cm]{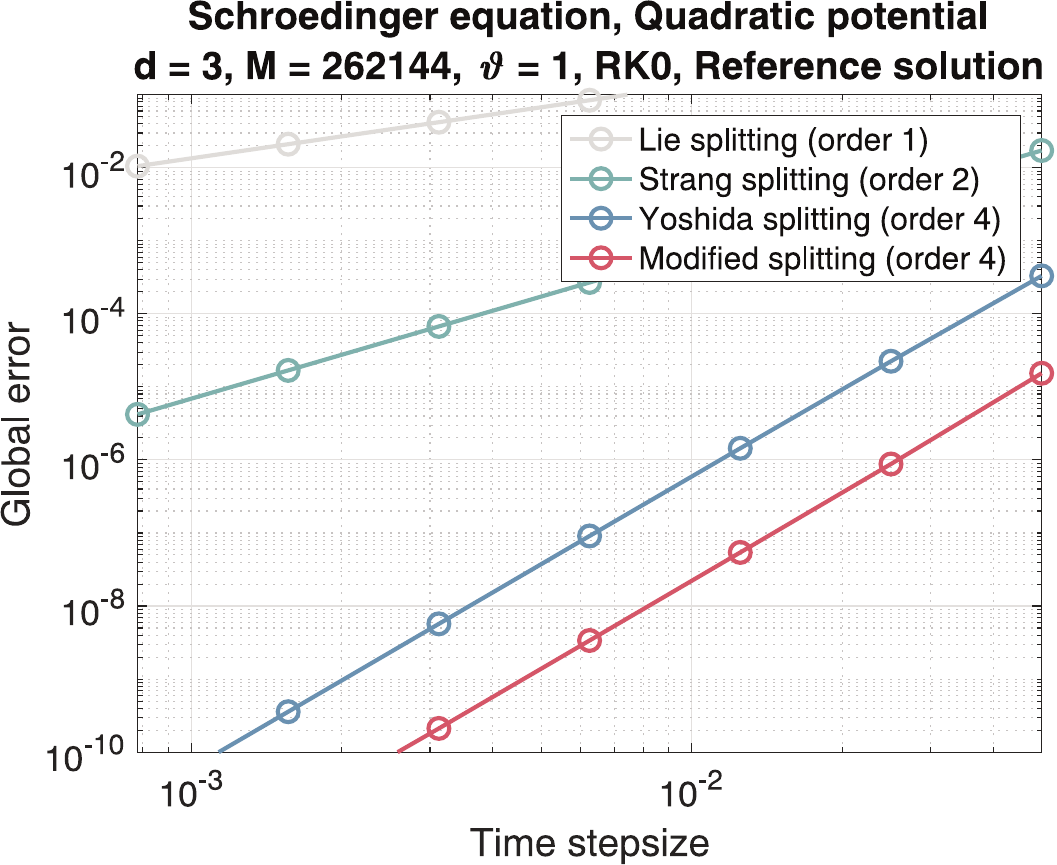} 
\end{center}
\caption{Time integration of the Gross--Pitaevskii equation~\eqref{eq:GPE} involving a quadratic potential by standard splitting methods and the novel modified operator splitting method.
Global errors versus time stepsizes in space dimensions $d \in \{1, 2, 3\}$.
Nonlinear ($\vartheta = 1$) versus simplified linear ($\vartheta = 0$) case.
Due to the validity of the invariance principle, the application of an explicit Runge--Kutta method is not needed (RK0).}
\label{fig:Figure1}
\end{figure}

\begin{figure}[t!]
\begin{center}
\includegraphics[width=6cm]{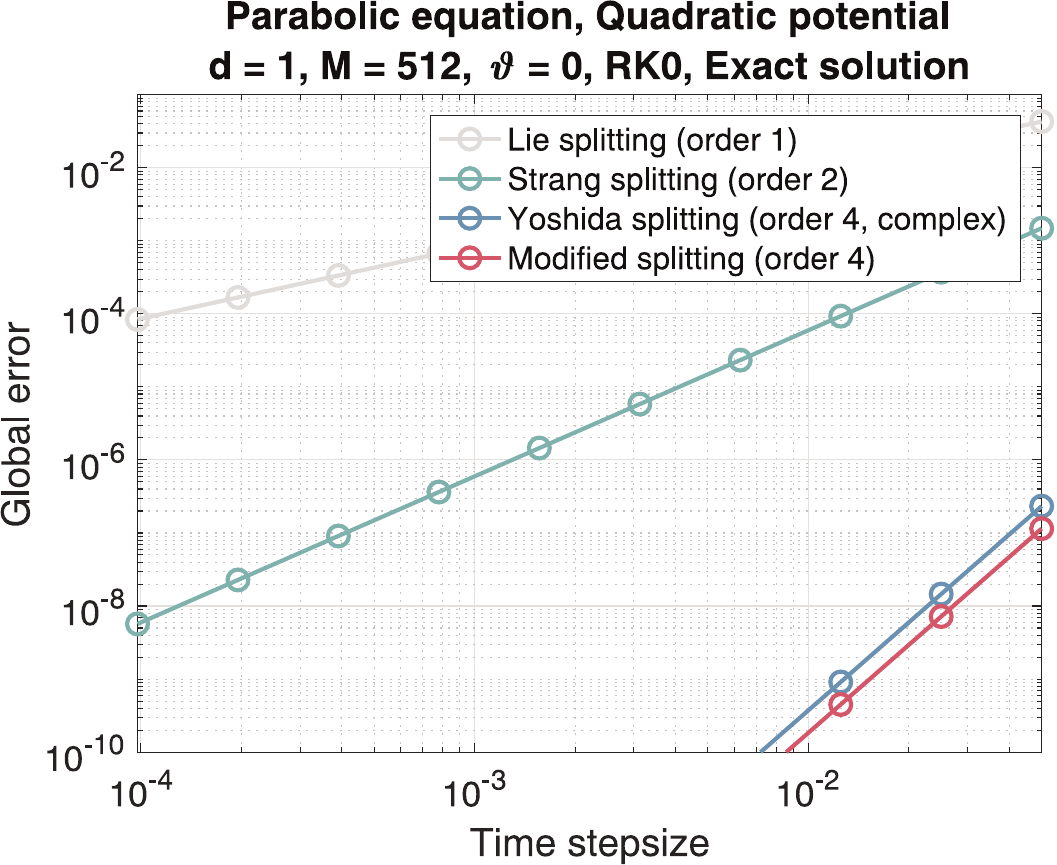} \quad
\includegraphics[width=6cm]{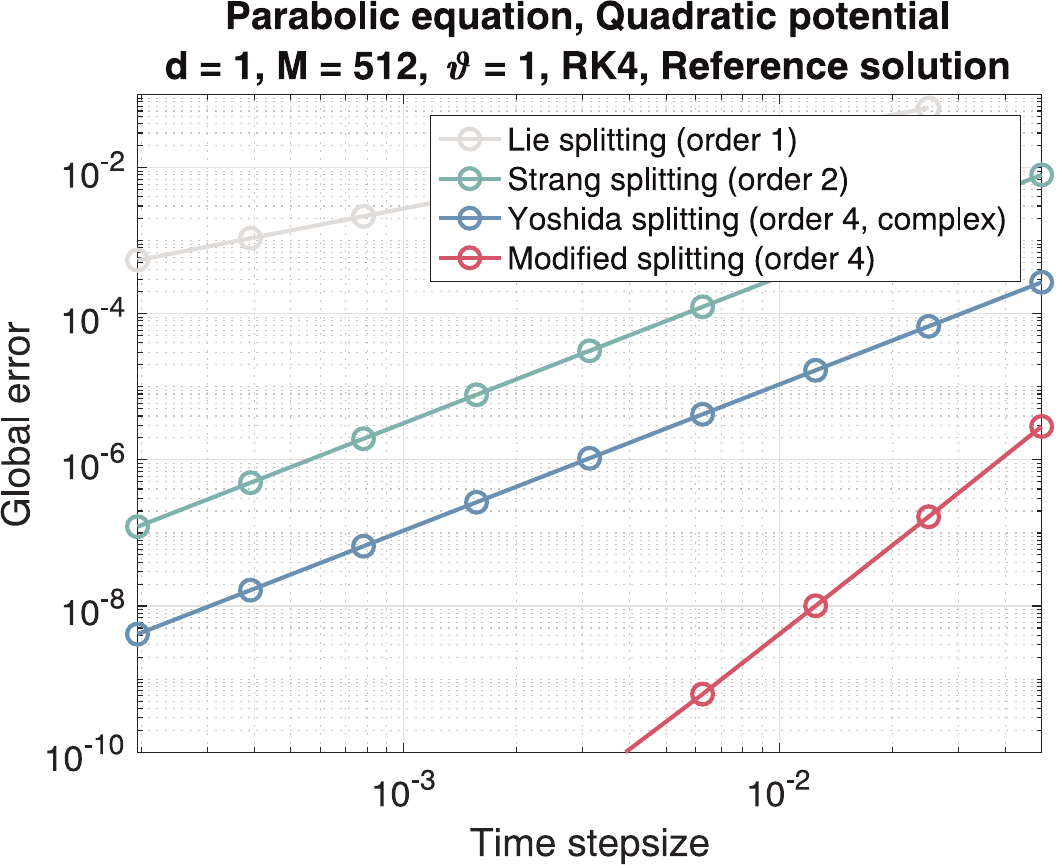} \\[4mm]
\includegraphics[width=6cm]{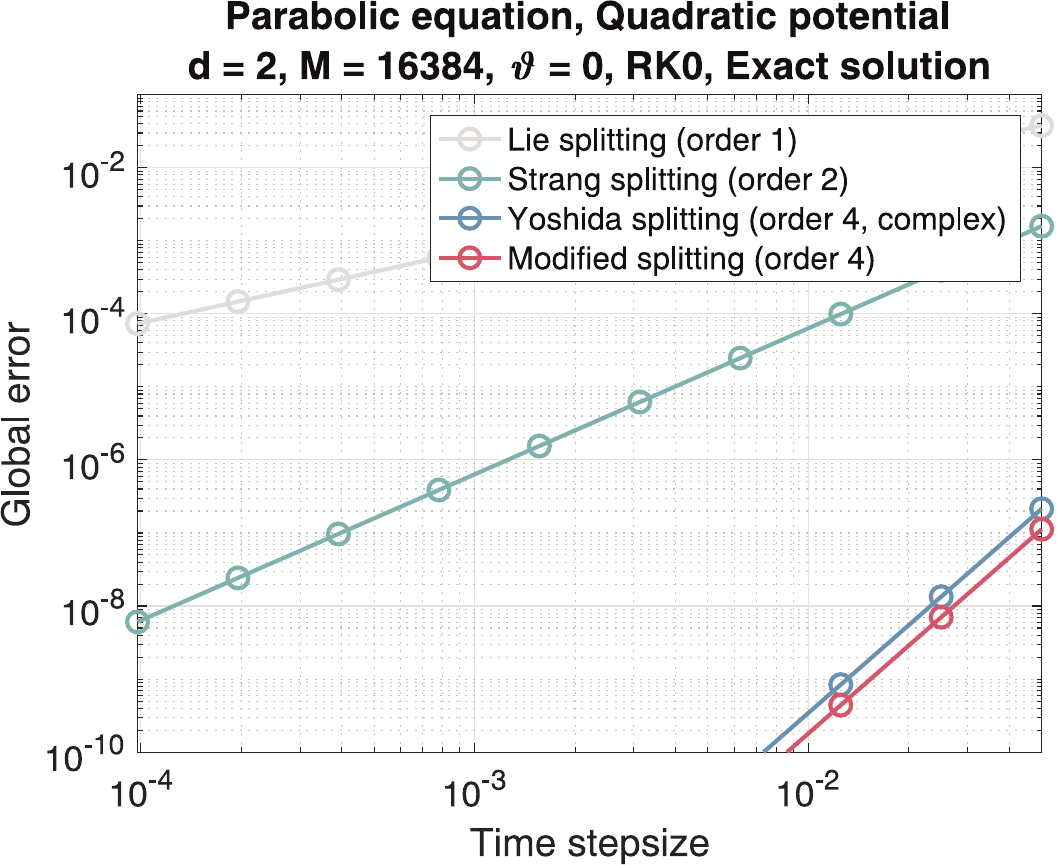} \quad
\includegraphics[width=6cm]{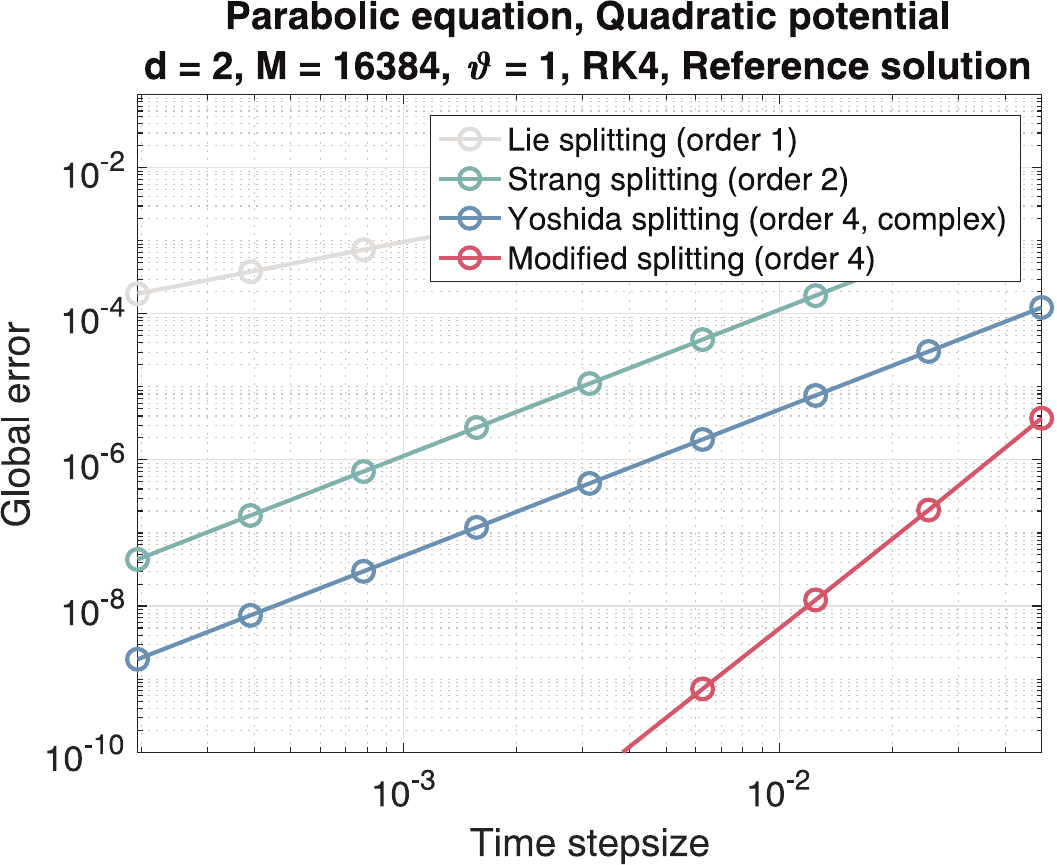} \\[4mm]
\includegraphics[width=6cm]{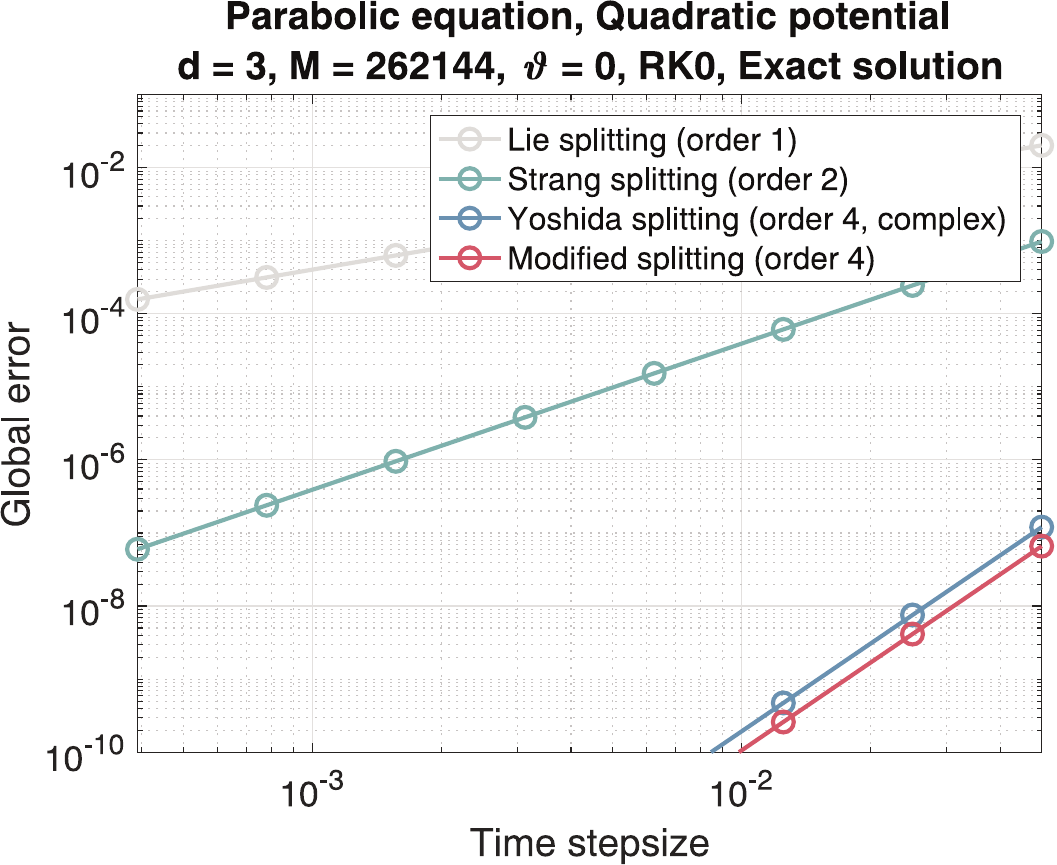} \quad
\includegraphics[width=6cm]{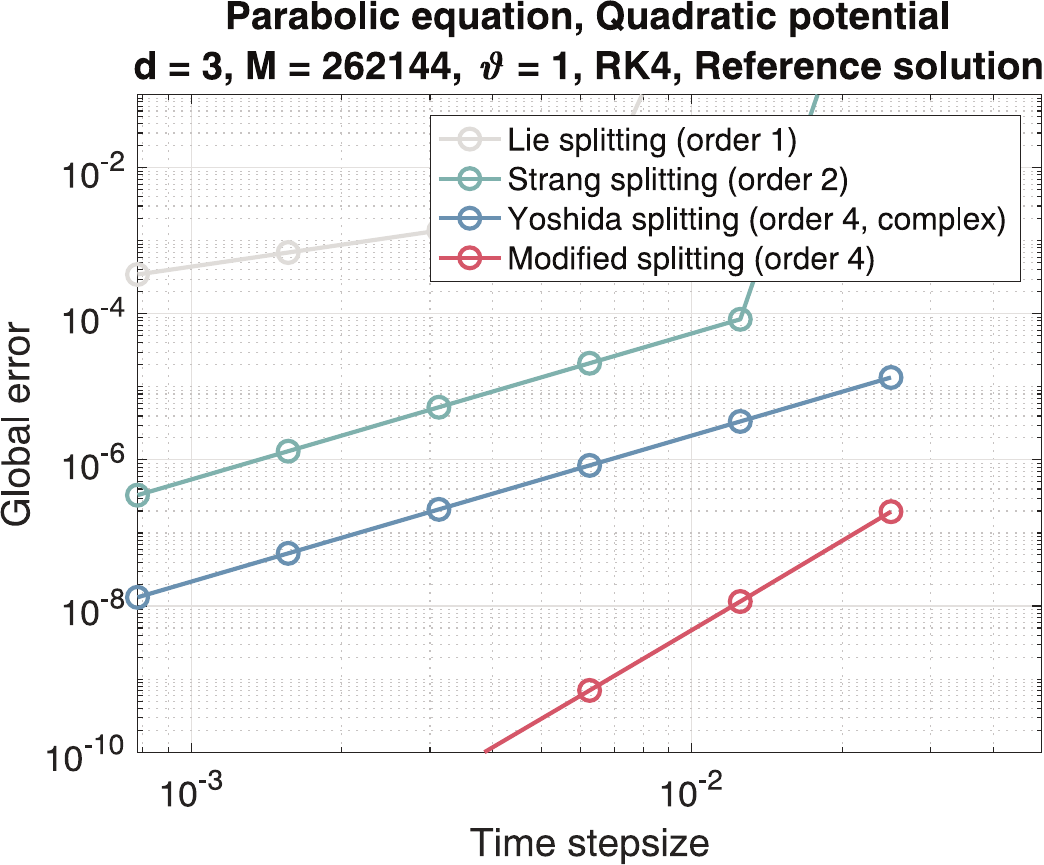} 
\end{center}
\caption{Time integration of the parabolic problem~\eqref{eq:Parabolic} involving a quadratic potential by standard splitting methods and the novel modified operator splitting method.
Global errors versus time stepsizes in space dimensions $d \in \{1, 2, 3\}$.
Nonlinear ($\vartheta = 1$) versus simplified linear ($\vartheta = 0$) case.
In order to resolve the nonlinear subproblem, a fourth-order explicit Runge--Kutta method is applied (RK4).
Depending on the stiffness of the equation, stability is ensured for sufficiently small time stepsizes. 
For a naive implementation of the Yoshida splitting method with complex coefficients, an order reduction is observed, see Appendix~\ref{sec:AppendixOrderReduction}.}
\label{fig:Figure2}
\end{figure}

\begin{figure}[t!]
\begin{center}
\includegraphics[width=6cm]{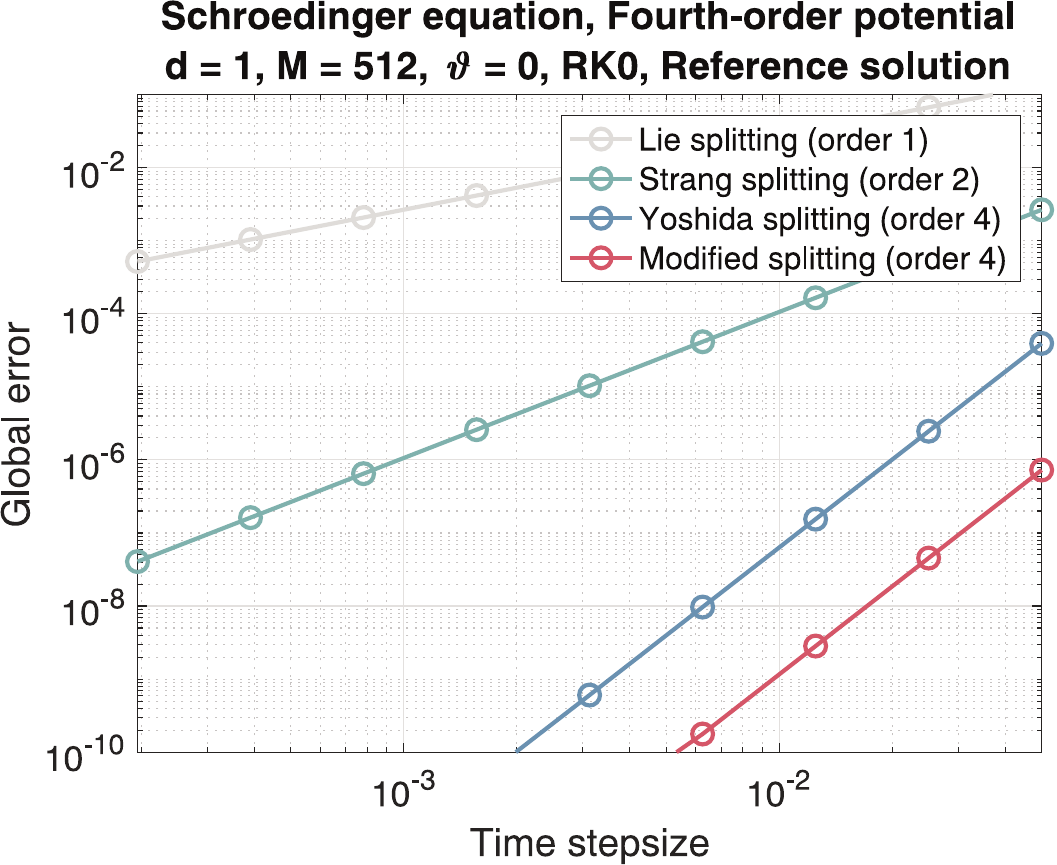} \quad
\includegraphics[width=6cm]{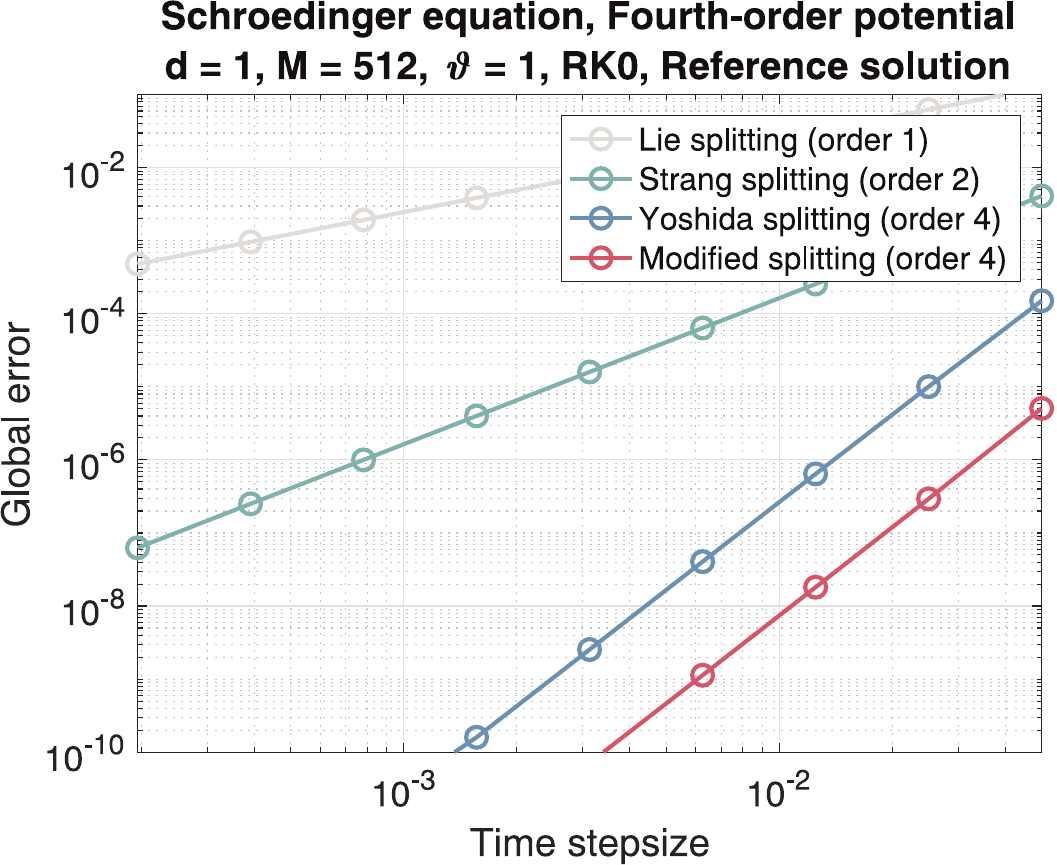} \\[4mm]
\includegraphics[width=6cm]{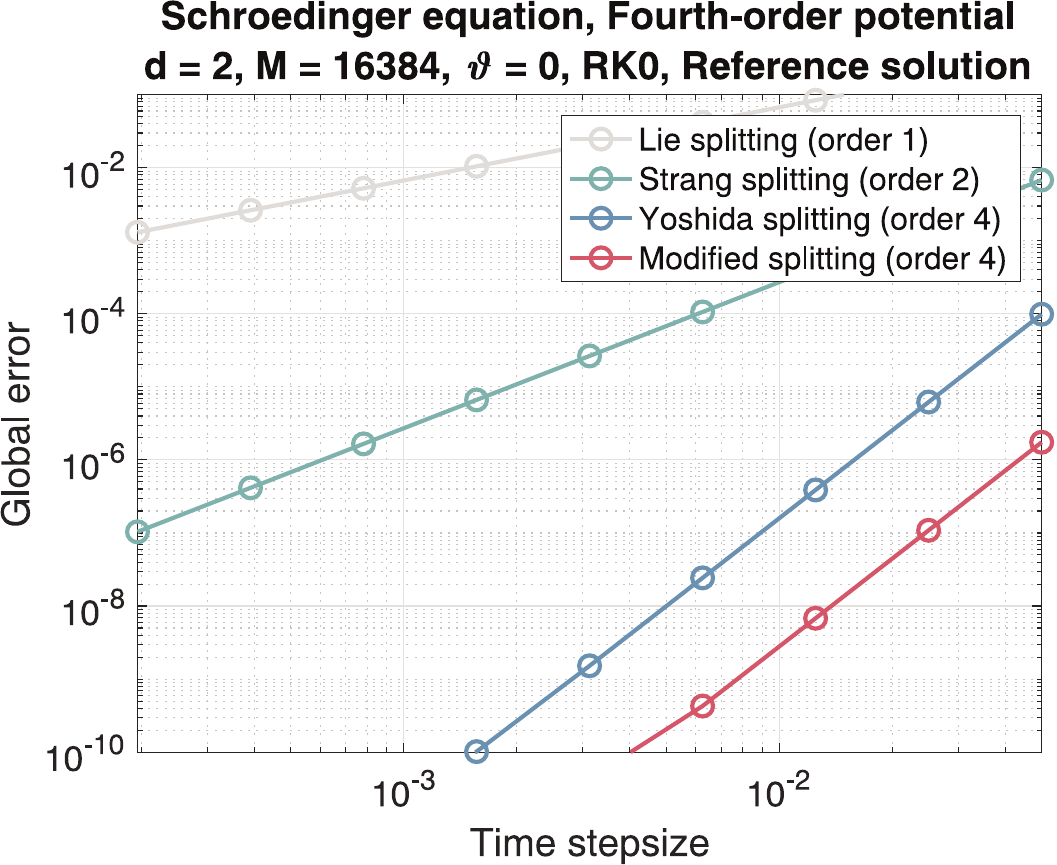} \quad
\includegraphics[width=6cm]{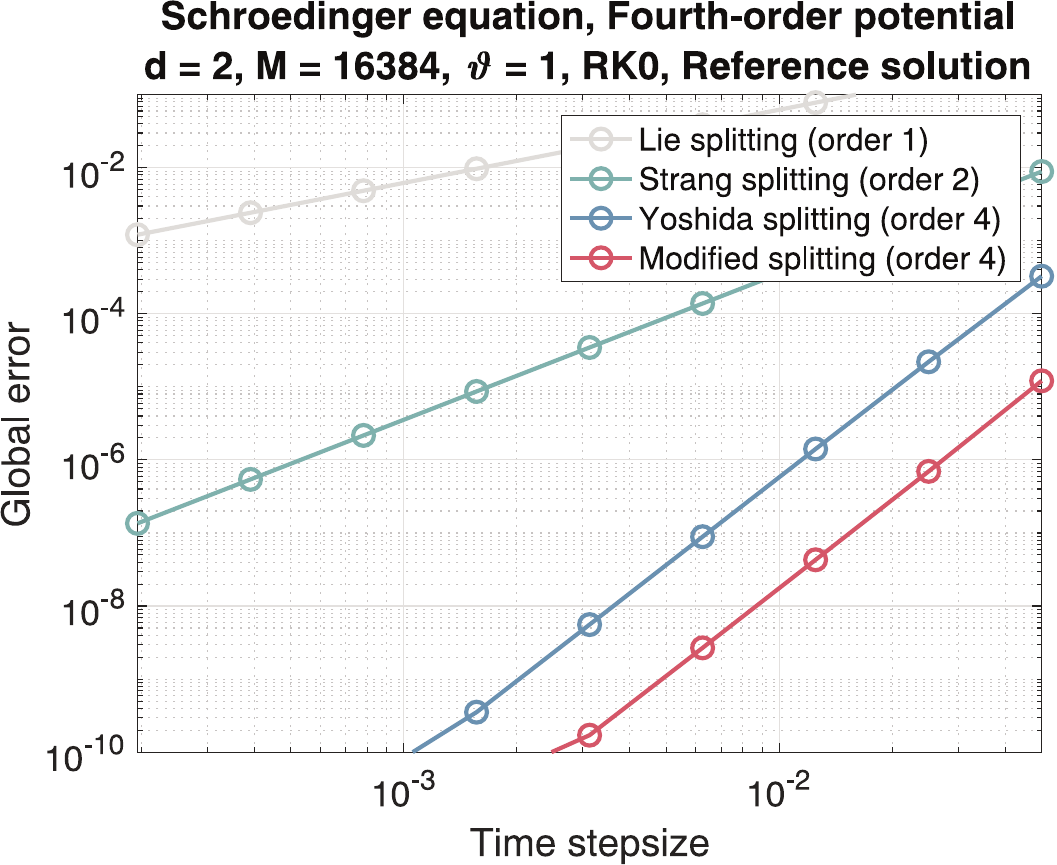} \\[4mm]
\includegraphics[width=6cm]{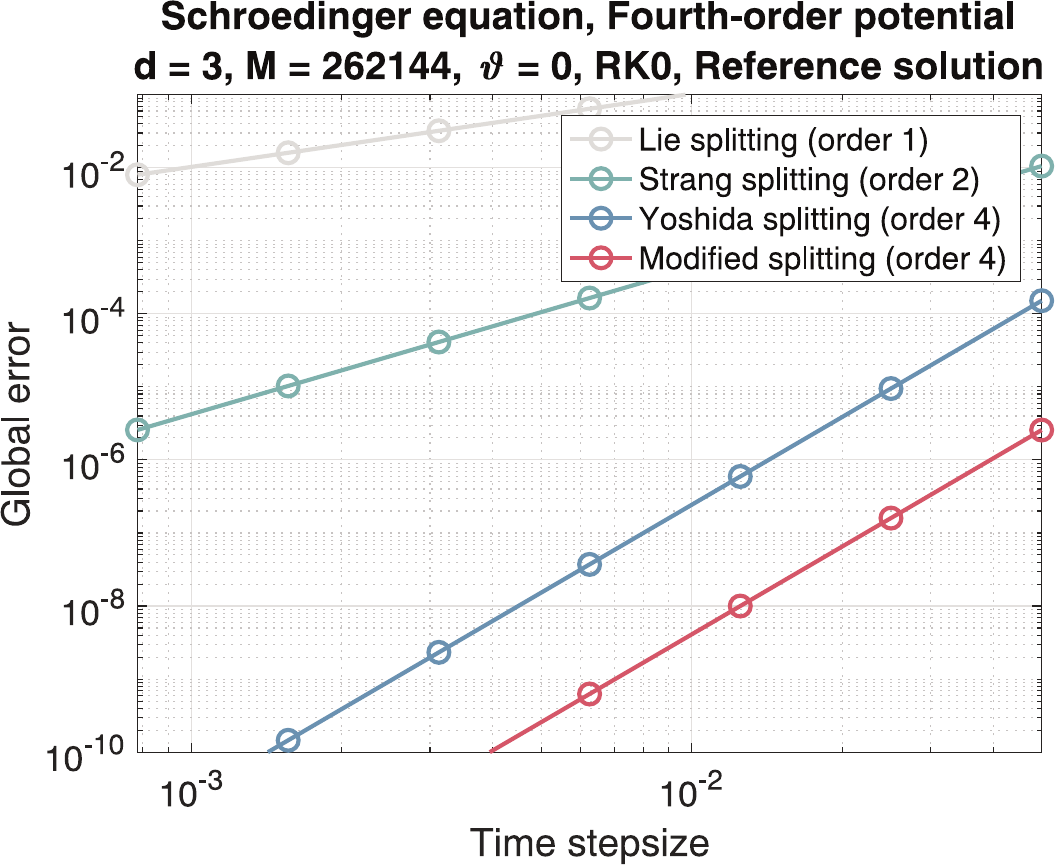} \quad
\includegraphics[width=6cm]{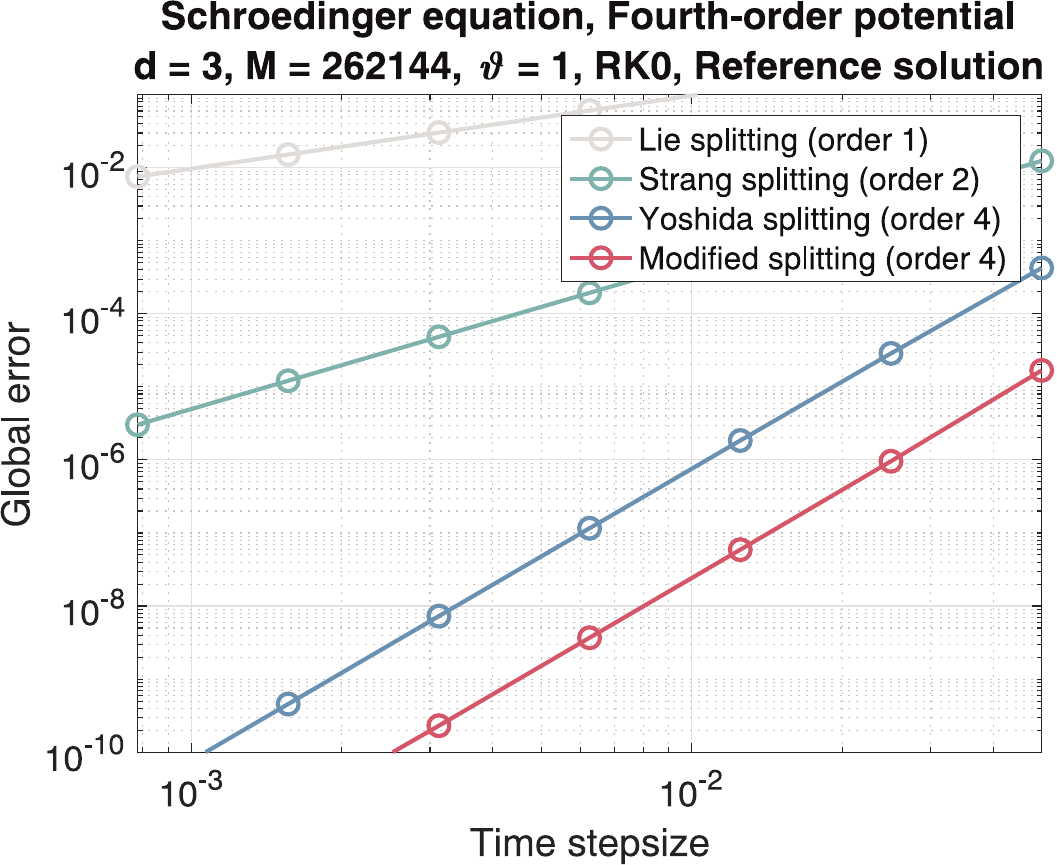} 
\end{center}
\caption{Corresponding results for the time-dependent Gross--Pitaevskii equation~\eqref{eq:GPE} involving a fourth-order polynomial potential.}
\label{fig:Figure3}
\end{figure}

\begin{figure}[t!]
\begin{center}
\includegraphics[width=6cm]{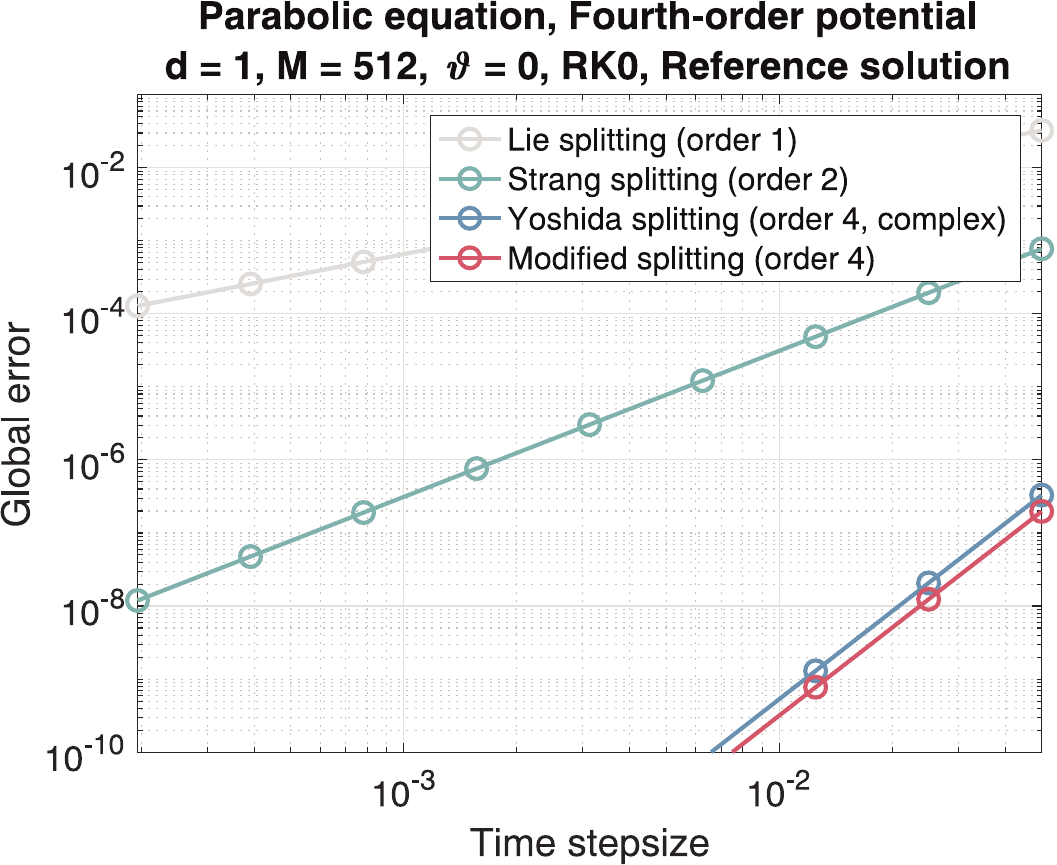} \quad
\includegraphics[width=6cm]{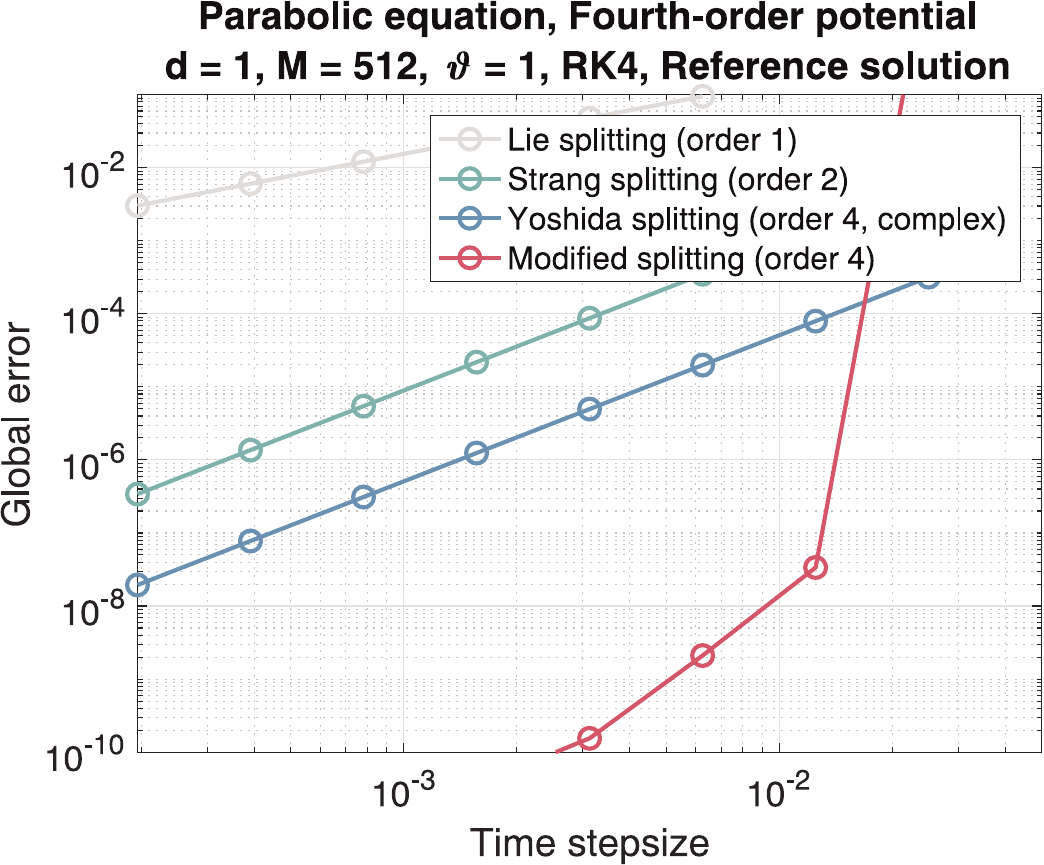} \\[4mm]
\includegraphics[width=6cm]{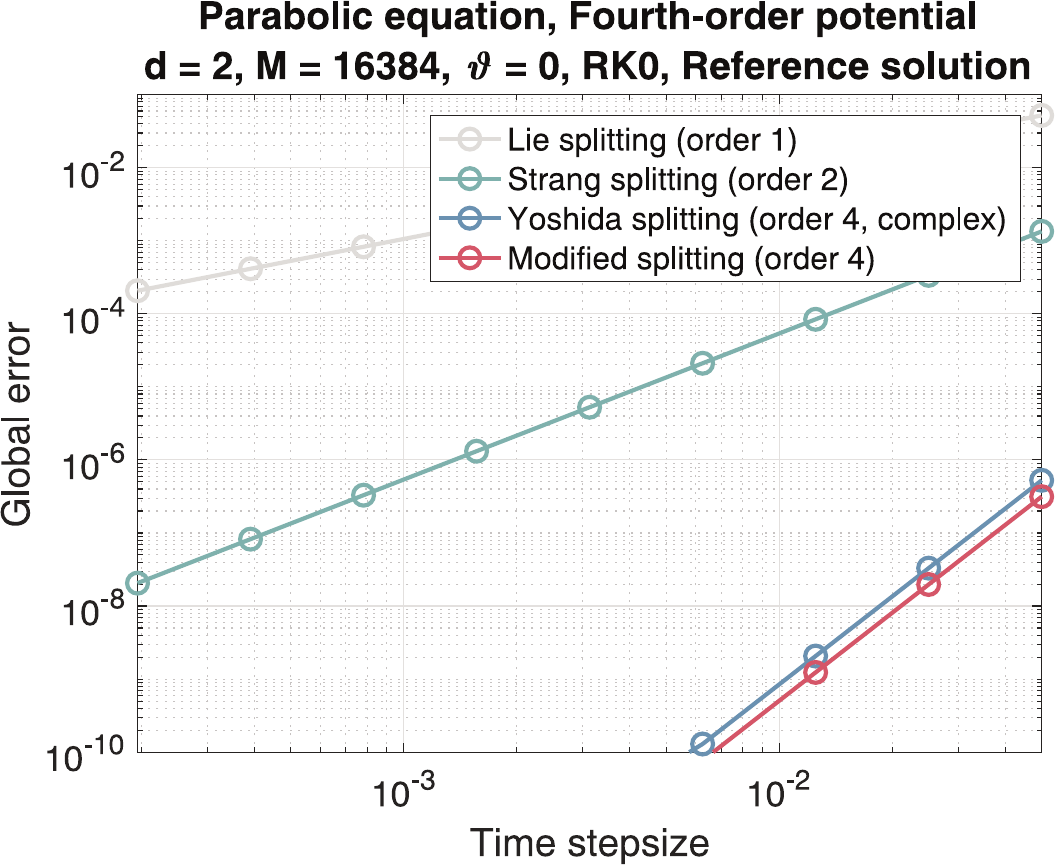} \quad
\includegraphics[width=6cm]{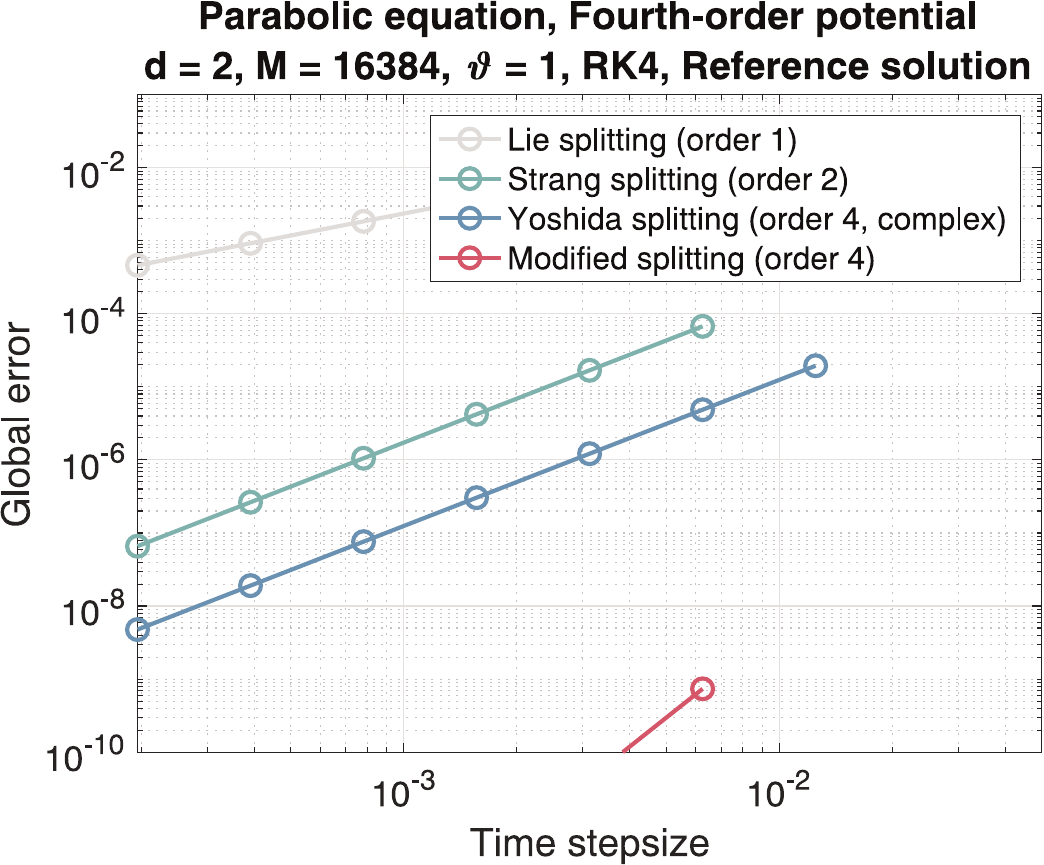} \\[4mm]
\includegraphics[width=6cm]{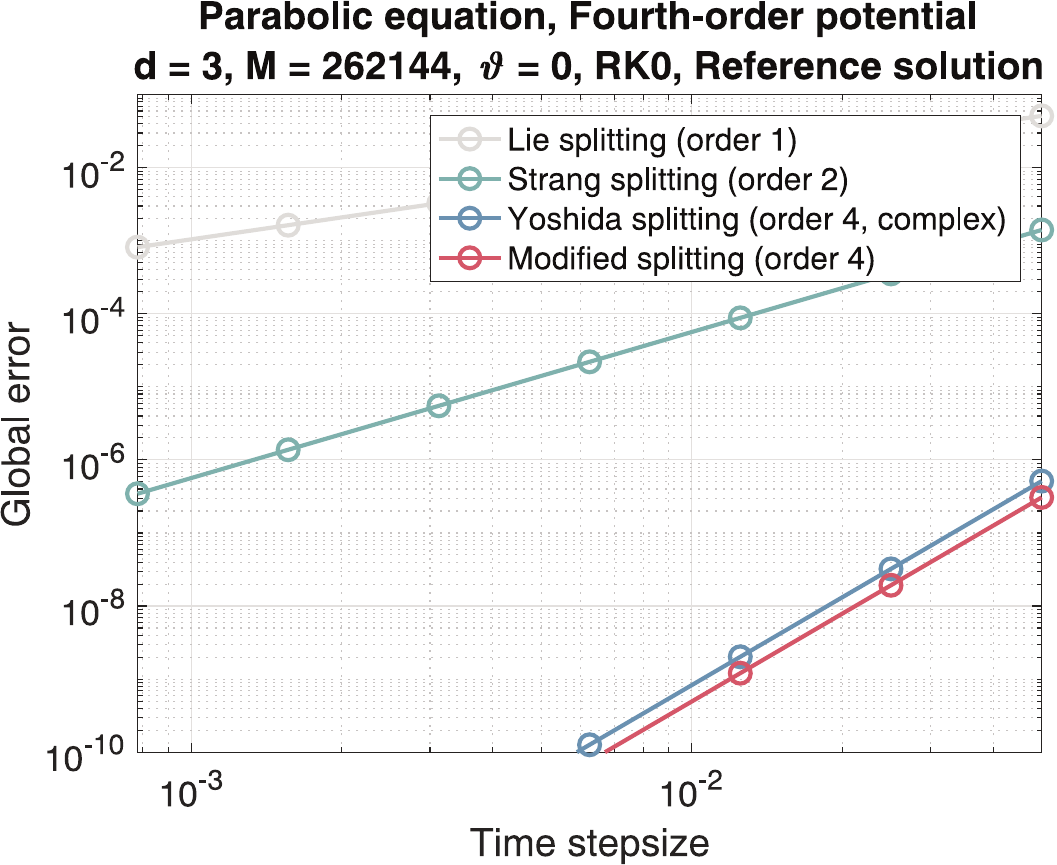}
\includegraphics[width=6cm]{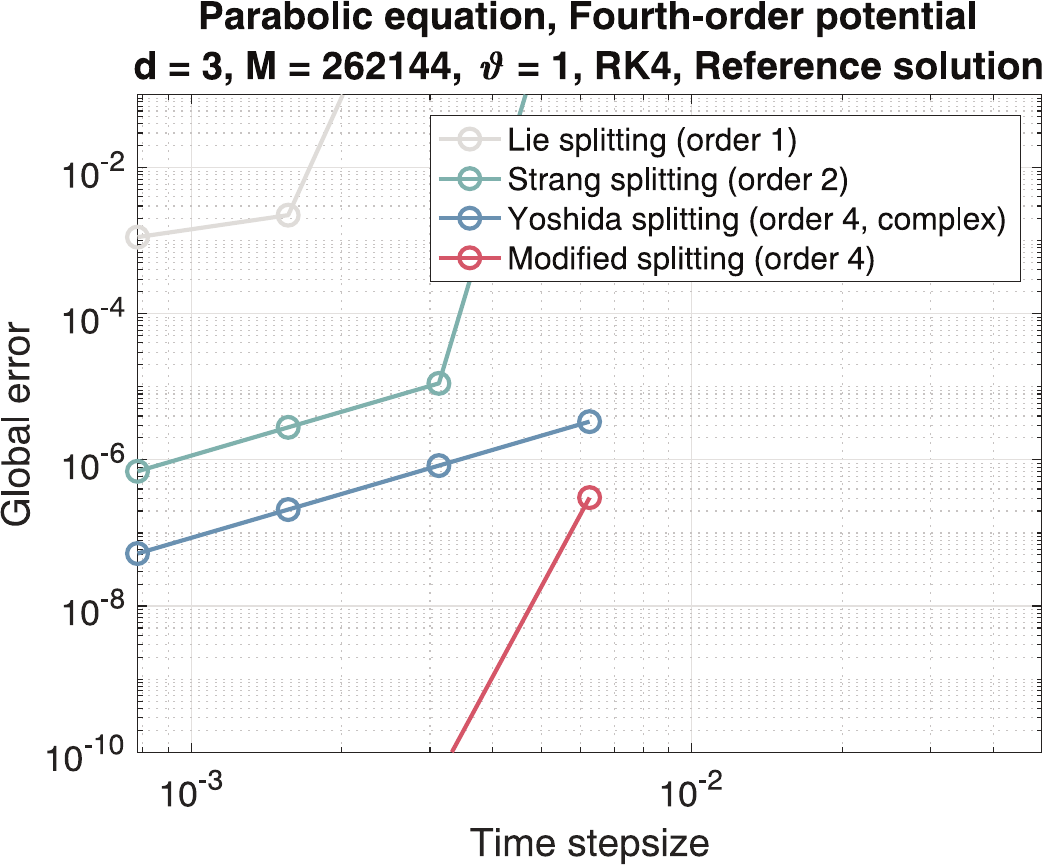} 
\end{center}
\caption{Corresponding results for the parabolic problem~\eqref{eq:Parabolic} involving a fourth-order polynomial potential.}
\label{fig:Figure4}
\end{figure}

\begin{figure}[t!]
\begin{center}
\includegraphics[width=6cm]{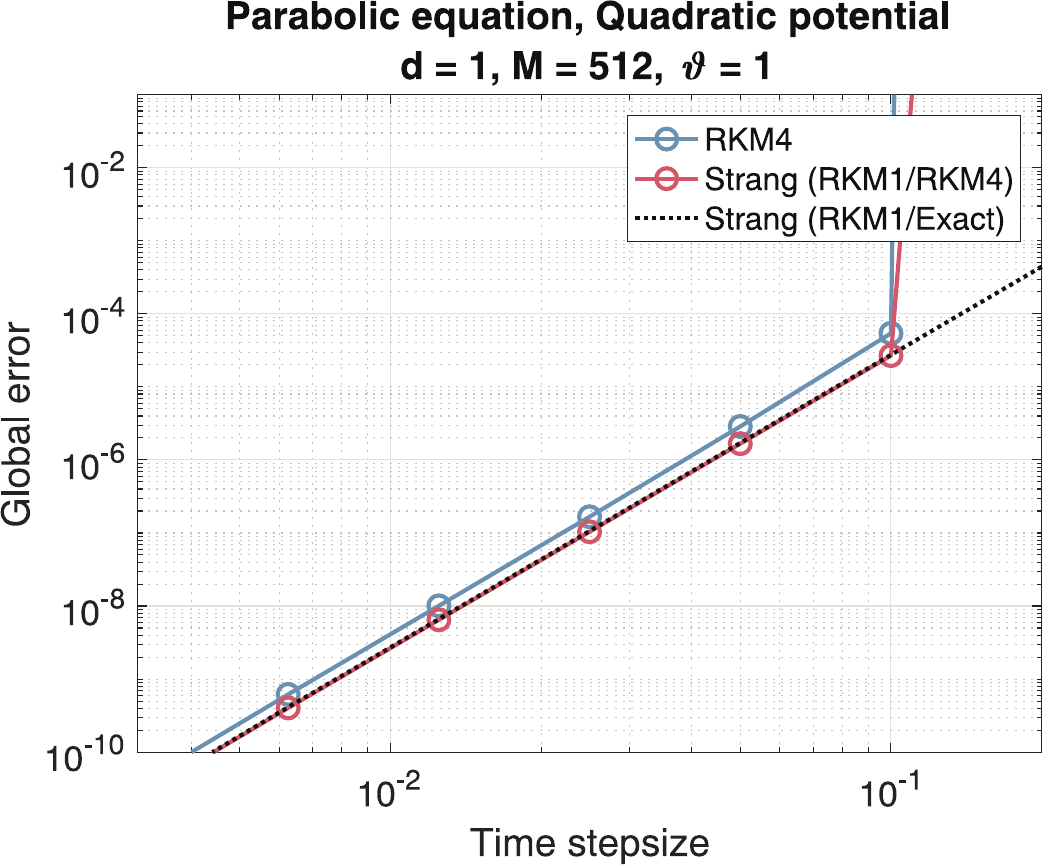} \quad
\includegraphics[width=6cm]{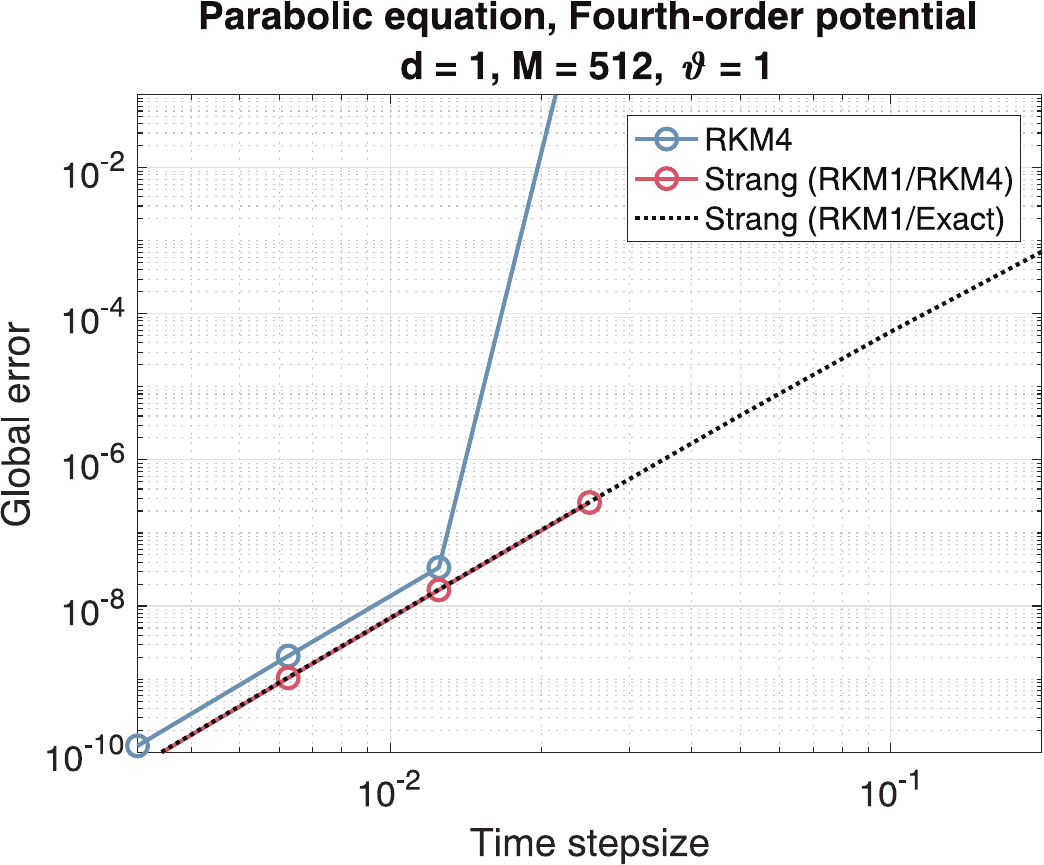}
\end{center}
\caption{Time integration of the one-dimensional parabolic equation~\eqref{eq:Parabolic} involving a quadratic potential (left) or a fourth-order potential (right), respectively, by the modified operator splitting method~\eqref{eq:ModifiedSplittingNonlinear}.
Global errors versus time stepsizes.
The original approach is based on the application of an explicit fourth-order Runge--Kutta method for the numerical solution of the nonlinear subproblem involving the double commutator
(cf. $\nE_{\tau, \frac{2}{3} F_2 - \frac{1}{72} \tau^2 G_2}$).
Alternative approaches are based on the Strang splitting method (cf. $\nE_{\frac{1}{2} \tau, \frac{2}{3} F_2} \circ \nE_{\tau, - \frac{1}{72} \tau^2 G_2} \circ \nE_{\frac{1}{2} \tau, \frac{2}{3} F_2}$).
Here, a reduced number of (inverse) fast Fourier transforms is required and an improved accuracy is observed.  
Furthermore, the knowledge of the exact solution to a component (cf. $\nE_{\frac{1}{2} \tau, \frac{2}{3} F_2} $) enhances the stability behaviour of the resulting time integration method for larger time increments.}
\label{fig:Figure5}
\end{figure}

\begin{figure}[t!]
\begin{center}
\includegraphics[width=10cm]{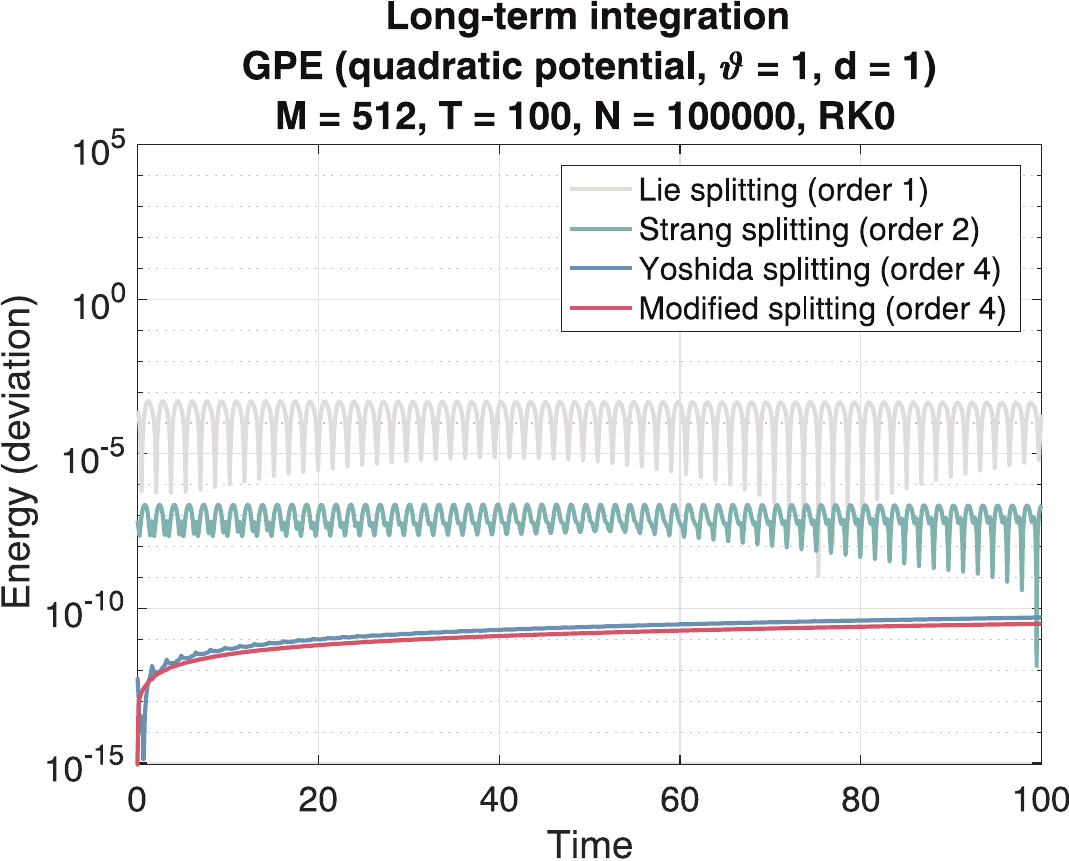} 
\end{center}
\caption{{Long-term integration of the one-dimensional Gross--Pitaevskii equation~\eqref{eq:GPE} by standard and modified operator splitting methods.
Computation of numerical approximations to the values of the energy at time grid points $t_n = n \, \tau$ for $\tau = 10^{-3}$ and $n \in \{0, 1, \dots, 10^5\}$ as well as corresponding deviations with respect to the minimal values, see also~\eqref{eq:Energy}.
The obtained results confirm the favourable geometric properties of the modified operator splitting method~\eqref{eq:ModifiedSplittingNonlinear}.}}
\label{fig:Figure6}
\end{figure}
\clearpage
\appendix 
\section{Matlab code}
\label{sec:AppendixMatlabCode}
A \textsc{Matlab} code that illustrates the practical implementation of modified operator splitting methods as outlined in Section~\ref{sec:Section6} and reproduces the numerical results displayed in Figures~\ref{fig:Figure1}--\ref{fig:Figure4} is available at 
\begin{center}
\url{doi.org/10.5281/zenodo.7945624}.
\end{center}
By default, the space dimenson is set equal to one such that the overall computation time amounts to a few minutes.
For dimensions two and three, respectively, the computational effort will increase accordingly to the complexity of the problem. 

\MyParagraph{Structure of the code}
The underlying process carries out the time integration of the Gross--Pitaevskii equation~\eqref{eq:GPE} and the related parabolic equation~\eqref{eq:Parabolic} by the standard splitting methods with coefficiens specified in~\eqref{eq:Coefficients} and the novel modified operator splitting methods~\eqref{eq:ModifiedSplittingNonlinear} with iterated commutators given in~\eqref{eq:G2Parabolic} and~\eqref{eq:G2GPE}, respectively.
The space discretisation is based on the Fourier spectral method and realised by fast Fourier transforms and their inverses. 
In loops over the selected space dimensions and the different test cases (parabolic / Schr{\"o}dinger equation, polynomial potential of degree two / four, linear case $\vartheta = 0$ / nonlinear case with $\vartheta = 1$), the global errors are computed for certain sequences of time stepsizes.
In order to make similarities apparent and contrast differences, several auxiliary functions are defined.
\begin{itemize}
\item 
\emph{Core}. \; 
Selection of time integration methods.
Definition of problem data (initial and final time, initial state, potential and its derivatives).
Choice of sequences of time stepsizes.
Performance of time integration.
Computation of global errors. 
\item 
\emph{PrecomputationFourier}. \;
Computation of underlying space grid and eigenvalues associated with Laplacian.
\item 
\emph{FourierReal2Spectral, FourierSpectral2Real}. \;
Fast Fourier transform and its inverse. 
\item 
\emph{PartA}. \;
Numerical solution of linear subproblems associated with Laplacian based on Fourier transforms.  
\item 
\emph{BWithoutU, B, TimeStepRKM124B, PartB}. \;
Numerical solution of nonlinear subproblems in the context of standard operator splitting methods.  
\item 
\emph{DoubleCommutator, BModifiedWithoutU, BModified, TimeStepRKM124\-BModified, PartBModified}. \;
Numerical solution of nonlinear subproblems in the context of modified operator splitting methods.  
\item 
\emph{TimeIntegration}. \;
Time integration by standard and modified operator splitting methods
\item 
\emph{MyTestCases}. \;
Definition of decisive quantities characterising different test cases.   
\item 
\emph{MyPlot}. \;
Visualisation of obtained results.
\end{itemize}

\MyParagraph{Specialisation and improvements}
The main purpose of the elementary structured \textsc{Matlab} code is the systematic comparison of the procedures required for the novel schemes with those for standard operator splitting methods, on the one hand
for Schr{\"o}dinger equations and on the other hand for parabolic equations.
We point out that a significant improvement of the performance will be achieved by separating the different types of evolution equations as well as time integrators and reconciling auxiliary functions.
Further enhancements in particular in connection with the more costly fast (inverse) Fourier transforms concern the distinction of real and complex arithmetics as well as the avoidance of redundancies, e.g., by reordering the eigenvalues associated with the Laplacian instead of the solution values and omitting additional scaling constants.  
\section{Order reduction of complex splitting methods}
\label{sec:AppendixOrderReduction}
In order to explain the observed order reduction for the complex splitting method~\eqref{eq:CoefficientsYoshidaComplex}, it suffices to study a nonlinear ordinary differential equation of the form~\eqref{eq:IVPNonlinear} with $F_1 = 0$ and $F_2(u) = \abs{u}^2 \, u$ on a single subinterval of length $(0, \tau)$.
More precisely, we consider the nonlinear subproblems 
\begin{equation*}
\begin{gathered}
\begin{cases}
\tfrac{\dd}{\dd t} \, u(t) = F_2^{(1)}\big(u(t)\big) = \big(u(t)\big)^3\,, \\
u(0) = u_0 \in \RR\,, \quad t \in (0, \tau)\,, 
\end{cases} 
\begin{cases}
\tfrac{\dd}{\dd t} \, u(t) = F_2^{(2)}\big(u(t)\big) = \abs{u(t)}^2 \, u(t)\,, \\
u(0) = u_0 \in \RR\,, \quad t \in (0, \tau)\,, 
\end{cases}
\end{gathered}
\end{equation*}
with coinciding real-valued solutions
\begin{equation*}
\nE_{\tau, F_2^{(1)}}(u_0) = u(\tau) = \nE_{\tau, F_2^{(2)}}(u_0)\,. 
\end{equation*}
On the one hand, using the Taylor series expansion
\begin{equation*}
u(\tau) = u(0) + \tau \, u'(0) + \tfrac{1}{2} \, \tau^2 \, u''(0) + \tfrac{1}{6} \, \tau^3 \, u'''(0) + \tfrac{1}{24} \, \tau^4 \, u''''(0) + \nO\big(\tau^5\big) 
\end{equation*}
in combination with the differential equation and derivatives thereof implies 
\begin{equation*}
u(\tau) = \Big(1 + \tau \, u_0^2 + \tfrac{3}{2} \, \tau^2 \, u_0^4 + \tfrac{5}{2} \, \tau^3 \, u_0^6 + \tfrac{35}{8} \, \tau^4 \, u_0^8\Big) \, u_0 + \nO\big(\tau^5\big)\,. 
\end{equation*}
As a consequence, for any real number $b \in \RR$, the relation 
\begin{equation*}
b \in \RR: \quad
u(b \, \tau) = \Big(1 + b \, \tau \, u_0^2 + \tfrac{3}{2} \, b^2 \, \tau^2 \, u_0^4 + \tfrac{5}{2} \, b^3 \, \tau^3 \, u_0^6 + \tfrac{35}{8} \, b^4 \, \tau^4 \, u_0^8\Big) \, u_0 + \nO\big(\tau^5\big)
\end{equation*}
is valid. 
For complex numbers, however, the more general expansion 
\begin{equation*}
\begin{split}
b, u_0 \in \CC: \quad
\nE_{\tau, b F_2^{(2)}}(u_0) 
&= \Big(1 + b \, \tau \, \abs{u_0}^2 + \big(b^2 + \tfrac{1}{2} \, \abs{b}^2\big) \, \tau^2 \, \abs{u_0}^4 \\
&\qquad + \big(b^3 + (\tfrac{7}{6} \, b + \tfrac{1}{3} \, \overline{b}) \, \abs{b}^2\big) \, \tau^3 \, \abs{u_0}^6 \\
&\qquad + \big(b^4 + \tfrac{1}{24} \, (46 \, b^2 + 6 \, \overline{b}^2 + 29 \, \abs{b}^2) \, \abs{b}^2\big) \, \tau^4 \, \abs{u_0}^8\Big) \, u_0 \\
&\qquad + \nO\big(\tau^5\big)
\end{split}
\end{equation*}
is obtained by decomposing the solution and accordingly the defining function into real and imaginary parts. 
In order to reproduce the approximation that corresponds to the Yoshida splitting, we impose the basic symmetry and order conditions
\begin{equation*}
b_4 = b_1\,, \quad b_3 = b_2\,, \quad b_1 + b_2 + b_3 + b_4 = 1\,,
\end{equation*}
and then perform the fourfold composition 
\begin{equation*}
u_{\text{Splitting}}(\tau) = \big(\nE_{\tau, b_1 F_2^{(2)}} \circ \nE_{\tau, b_2 F_2^{(2)}} \circ \nE_{\tau, b_2 F_2^{(2)}} \circ \nE_{\tau, b_1 F_2^{(2)}}\big)(u_0)\,. 
\end{equation*}
Requiring this expansion to coincide with the expansion of the exact solution
\begin{equation*}
u_{\text{Splitting}}(\tau) - u(\tau) = C(u_0) \, \Im(b_1) \, \big(1 - \Re(b_1)\big) \, \tau^3 + \nO\big(\tau^4\big)
\end{equation*}
leads to a condition that obviously contradicts the order conditions for splitting methods and explains the observation of local order three and global order two.
Similar arguments apply to evolution equations of Schr{\"o}dinger type such as the Gross--Pitaevskii equation.
But, in this context, the invariance principle permits to avoid the application of a fourth-order Runge--Kutta method to the arising nonlinear subproblem
\begin{equation*}
\tfrac{\dd}{\dd t} \, u(t) = F_2^{(2)}\big(u(t)\big) = - \, \ii \, \abs{u(t)}^2 \, u(t)\,, \quad t \in (0, \tau)\,.
\end{equation*}
\end{document}